\pdfoutput=1
\documentclass[%
11pt,
onecolumn,
tightenlines,
superscriptaddress,
preprintnumbers,
nofootinbib,
amsmath,amssymb,amsthm,
physrev,
eqsecnum,tikz,
]{revtex4-2}

\usepackage{isomath}
\usepackage{amsmath,amsthm}
\usepackage{amsbsy}
\usepackage{amssymb}
\usepackage{amscd}
\usepackage{amsfonts}
\usepackage{stmaryrd}
\usepackage{siunitx}
\usepackage{euscript}
\usepackage[utf8]{inputenc}
\usepackage[T1]{fontenc}
\usepackage{newtxtext} 
\everymath{\displaystyle}
\usepackage{exscale}

\usepackage{graphicx}
\usepackage{boxedminipage}
\usepackage{calc}
\usepackage[usenames,dvipsnames]{xcolor}
\graphicspath{ {media/} }
\usepackage[caption=false,justification=centerlast]{subfig}

\usepackage{setspace}
\usepackage{enumitem}
\setitemize{noitemsep,topsep=0pt,parsep=0pt,partopsep=0pt}
\setenumerate{noitemsep,topsep=0pt,parsep=0pt,partopsep=0pt}
\setdescription{noitemsep,topsep=0pt,parsep=0pt,partopsep=0pt}

\usepackage{soul} 

\usepackage{hyperref}
\usepackage[normalem]{ulem}

\usepackage[small]{titlesec}

\titlespacing*{\section}{0pt}{12pt plus 4pt minus 2pt}{2pt plus 2pt minus 2pt}
\titlespacing*{\subsection}{0pt}{12pt plus 4pt minus 2pt}{2pt plus 2pt minus 2pt}
\titlespacing*\subsubsection{0pt}{12pt plus 4pt minus 2pt}{2pt plus 2pt minus 2pt}
\titlespacing*\paragraph{0pt}{12pt plus 4pt minus 2pt}{2pt plus 2pt minus 2pt}

\makeatletter

    \renewcommand*{\p@subsection}{}
    
    \renewcommand*{\p@subsubsection}{}
\makeatother

\usepackage{isomath}
\usepackage{amsmath}
\usepackage{amssymb}
\usepackage{amscd}
\usepackage{amsfonts}

\newcommand{\bfsigma}{\mathbold {\sigma}}
\newcommand{\bfepsilon}{\mathbold {\epsilon}}

\DeclareMathOperator{\divergence}{div}

\DeclareMathOperator{\trace}{tr}

\newcommand{\parderiv}[2]{\frac{\partial #1}{\partial #2}}
\newcommand{\dm}{\ \mathrm{d}}
\newcommand{\deriv}[2]{\frac{\dm #1}{\dm #2}}

\newcommand{\bfb}{{\mathbold b}}

\newcommand{\bfg}{{\mathbold g}}

\newcommand{\bfk}{{\mathbold k}}

\newcommand{\bfn}{{\mathbold n}}

\newcommand{\bfq}{{\mathbold q}}

\newcommand{\bft}{{\mathbold t}}
\newcommand{\bfu}{{\mathbold u}}
\newcommand{\bfv}{{\mathbold v}}

\newcommand{\bfx}{{\mathbold x}}

\newcommand{\bfF}{{\mathbold F}}

\newcommand{\bfI}{{\mathbold I}}

\newcommand{\bfK}{{\mathbold K}}

\newcommand{\bfT}{{\mathbold T}}

\begin{document}

\preprint{To appear in International Journal for Numerical and Analytical Methods in Geomechanics}
\preprint {DOI: \url{ https://doi.org/10.1002/nag.3326}}

\title{\Large{Energetic Formulation of Large-Deformation Poroelasticity}}

\author{Mina Karimi}
    \email{minakari@andrew.cmu.edu}
    \affiliation{Department of Civil and Environmental Engineering, Carnegie Mellon University}

\author{Mehrdad Massoudi}%
    \affiliation{National Energy Technology Laboratory, Pittsburgh PA 15236-0940}%

\author{Noel Walkington}
    \affiliation{Center for Nonlinear Analysis, Department of Mathematical Sciences, Carnegie Mellon University}

\author{Matteo Pozzi}
    \affiliation{Department of Civil and Environmental Engineering, Carnegie Mellon University}
    
\author{Kaushik Dayal}
    \email{Kaushik.Dayal@cmu.edu}
    \affiliation{Department of Civil and Environmental Engineering, Carnegie Mellon University}
    \affiliation{Center for Nonlinear Analysis, Department of Mathematical Sciences, Carnegie Mellon University}
    \affiliation{Department of Materials Science and Engineering, Carnegie Mellon University}
    
\date{\today}


\begin{abstract}
    The modeling of coupled fluid transport and deformation in a porous medium is essential to predict the various geomechanical process such as CO2 sequestration, hydraulic fracturing, and so on.
    Current applications of interest, for instance, that include fracturing or damage of the solid phase, require a nonlinear description of the large deformations that can occur.
    This paper presents a variational energy-based continuum mechanics framework to model large-deformation poroelasticity.
    The approach begins from the total free energy density that is additively composed of the free energy of the components.
    A variational procedure then provides the balance of momentum, fluid transport balance, and pressure relations.
    A numerical approach based on finite elements is applied to analyze the behavior of saturated and unsaturated porous media using a nonlinear constitutive model for the solid skeleton.
    Examples studied include the Terzaghi and Mandel problems; a gas-liquid phase-changing fluid; multiple immiscible gases; and unsaturated systems where we model injection of fluid into soil.
    The proposed variational approach can potentially have advantages for numerical methods as well as for combining with data-driven models in a Bayesian framework.
\end{abstract}

\maketitle



\section{Introduction} \label{Introduction}

The modeling of fluid transport in porous deformable media, and analyzing the coupled deformation-transport behavior, is of importance for a variety of engineering and scientific applications, ranging from the modeling of natural and engineered biological materials to geomechanics, e.g. \cite{ehlers2009challenges,song2020spatially,jin2020fluid,selvadurai2013mechanics,simon1992multiphase,duddu2020non,alaimo2015laminar,sun2013stabilized,taleghani2020thermoporoelastic,lu2018modeling,fan2019basement,duda2010theory,von2020morphogenesis,mudunuru2012framework,srinivasan2019model,trivisa2018analysis,fritz2021modeling,jha2014coupled,simon1992multiphase,hassanizadeh1979general,bui2020scalable} and numerous other works.
It has therefore been the focus of numerous formulations, starting from the seminal work of Biot \cite{biot1935problem,biot1941general}.
In this paper, we formulate a general finite-deformation model for poroelasticity based on a variational approach.

The classical linear theory proposed by Biot has been applied to study a wide range problems, e.g. \cite{coussy2004poromechanics,wang2000theory,rudnicki2001coupled,cowin1999bone}.
However, the many assumptions of the linear theory, including linear elastic response for the solid phase, infinitesimal strain, and constant permeability are too restrictive to model the complex nonlinear response of many porous systems, such as polymers, damaged geomaterials with fracture, and soft tissues. 

Therefore, there have been many efforts to generalize the theory to the nonlinear setting, with early works starting from Biot \cite{biot1973nonlinear}.
More recently, examples of nonlinear large-deformation poroelastic models include models for elastomeric and biological materials \cite{macminn2016large}; multiphase continuum formulations presented in \cite{federico2012elasticity} and \cite{hong2008theory} for saturated fields reinforced by elastic fibers and polymeric gels, respectively; 
and for computational geomechanics, such as multiphase finite element formulations for fully and partially saturated porous geomaterials, by Borja and coworkers \cite{li2004dynamics,borja2004cam,song2014mathematical}.
Another older approach is based on mixture theory, e.g. \cite{bowen1982compressible}, which leads to a formidable set of equations.

The key contribution of this paper is a variational derivation that starts from a free energy, and uses variational principles to obtain the governing equations, in contrast to other works in the literature.
For instance, a closely-related body of work is by Borja and coworkers \cite{borja2006mechanical, borja2004cam}, in which they use a thermodynamic approach based on an energy and also the balance laws of continuum mechanics; they formulate a three-phase continuum mixture theory to analyze a partially saturated porous media. 
They use a free energy function of the medium is defined based on the solid strain, displacement vector, and pressure of each phase. 
In a related approach, \cite{kleinfelter2007mixture} formulated an energy-based method for unsaturated porous medium where they used the energy equation which is developed by \cite{bennethum1996multiscale}. 

Our approach is instead to directly use the energy and associated variational principles to derive the governing equations.
Similar approaches have been used in geomechanics to model solid materials with a fluid phase, i.e. two-phase porous medium models, e.g. \cite{borja1995mathematical,semnani2016thermoplasticity,anand20152014}.
Variational approaches for poromechanics are relatively recent, and appear to begin with \cite{gajo2010general,chester2010coupled}.
While we follow the overall structure of these works in our formulation, these have been restricted to simpler settings, e.g. saturated systems, two-phase systems, and so on.
In this work, as in those works, the governing equations are derived from energy minimization.
However, \cite{chester2010coupled} assume a multiplicative decomposition of the deformation gradient into an elastic part and a poromechanical part -- corresponding to the deformation caused by fluid transport -- whereas we do not use this assumption.

This work presents a variational energy-based model for poromechanics, which can be applied to both saturated and unsaturated porous mediums, containing compressible constituents. 
The proposed model is amenable to the use of arbitrary energy density functions for both the fluid and solid phases. 
In particular, it can model unsaturated porous media with compressible constituents, and multiple immiscible fluid phases with  compressible and incompressible constituents.
Our strategy is to begin with the free energy density function, and use a variational approach to obtain the conservation of mass and momentum, and the pressure equality equations.
An important advantage of the energetic formulation and the minimization structure is the advantage for robust numerical schemes using finite elements and other methods \cite{strang}.
In addition, it can enable efficient numerical approaches to combining physics-based models with data-driven strategies \cite{karimi2022}.

\paragraph*{Organization.}
Section \ref{sec:Model Formulation} provides the general formulation, starting with the kinematics and then using a variational approach to develop the general governing equations.
Sections \ref{Constitutive Models for Saturated Systems} discusses the constitutive models for saturated and unsaturated systems, which includes a discussion on porous media with single and multiple fluids.
Section \ref{Models for Referential and Current Volume Fractions} examines the relation between the referential and current volume fractions for different assumptions on the fluid and solid responses.
Section \ref{Model verification} provides numerical solutions of benchmark problems proposed by Terzaghi and Mandel. 
Sections \ref{single compressible fluid}, \ref{multiple compressible fluid}, and \ref{unsaturated system with incompressible fluid} provides numerical examples to investigate the settings of a phase-changing (gas to liquid) fluid, multiple compressible fluids, and unsaturated systems.


\section{Model Formulation} \label{sec:Model Formulation}

We assume that our system consists of $N$ phases or components indexed by $i$, with $i=1$ corresponding to the solid skeleton phase, and $i = 2, \ldots, N$ corresponding to the fluid phases.
We follow the motion of the solid skeleton with respect to the reference configuration, i.e., we use a Lagrangian description. Therefore, we use material time derivatives. Once we obtain the governing equations, they can be pushed forward to the current configuration.

The gradient and divergence operators in the reference and current configurations will be denoted $\nabla_0, \nabla, \divergence_0, \divergence$ respectively.
In general, quantities with the subscript $0$ will refer to referential quantities.
We use $\deriv{\ }{t}$ to denote the material time derivative following the motion of the solid skeleton.

\subsection{Kinematics} \label{sec:Kinematics}

\begin{figure}[htb!]
		\includegraphics[width=0.8\textwidth]{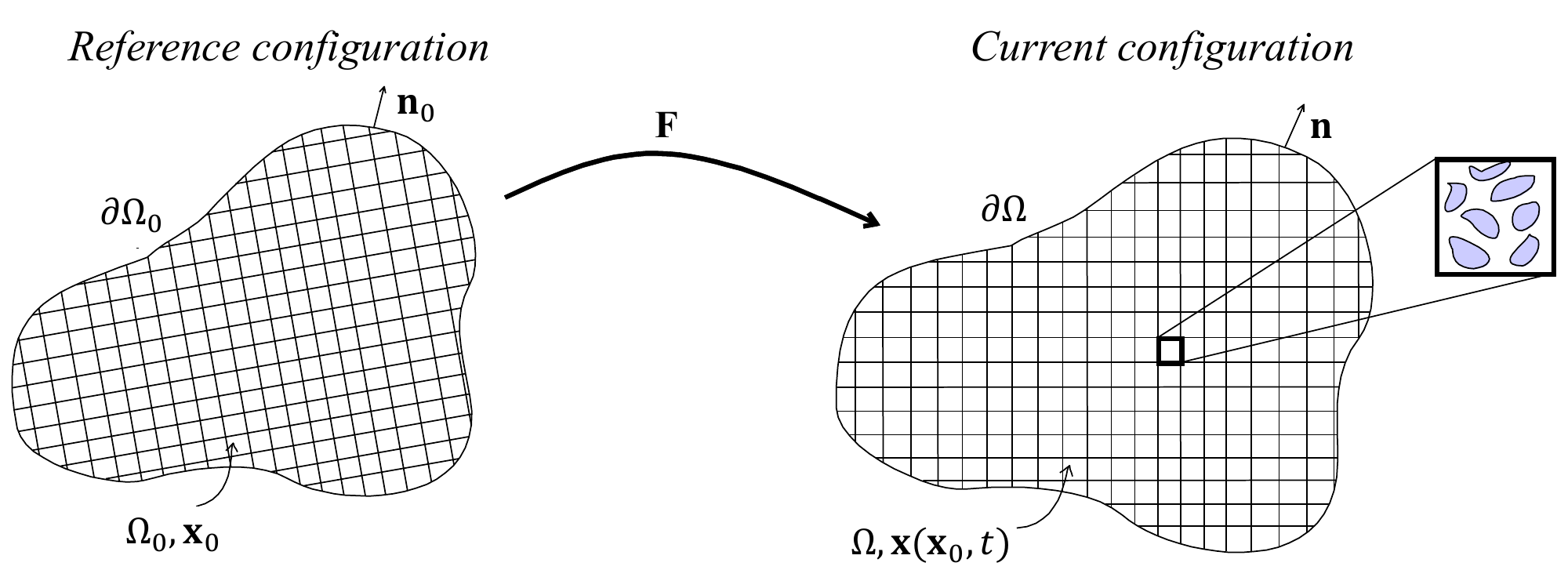}
		\caption{Representative volume element of the porous medium, and relation of reference and current configurations.}
		\label{fig:2-scale}
\end{figure}

We consider a porous medium $\Omega$, containing solid and fluid components.
We consider a two-scale picture motivated by homogenization.
Namely, there is a macroscopic deformation $\bfx(\bfx_0,t)$, where $\bfx$ is the location at time $t$ of a material particle located at $\bfx_0$ in the reference configuration.
To each macroscale location, i.e., for a given $\bfx_0$, there corresponds a microscale representative volume element (RVE), which accounts for the geometry of the porous skeleton and other microscale features (Figure \ref{fig:2-scale}).
Therefore, at the macroscopic homogenized, an RVE is equivalent to a material point.

We define:
\begin{enumerate}
    \item the macroscale deformation gradient tensor $\bfF:=\parderiv{\bfx}{\bfx_0}$, the Jacobian $J:=\det\bfF$, and the displacement $\bfu(\bfx_0,t):=\bfx(\bfx_0,t) - \bfx_0$;
    \item the densities $P_i$ and $P_{0i}$ as the mass of component $i$ per unit (current and referential respectively) volume of the mixture;
    \item the true densities $\rho_i$ and $\rho_{0i}$ as the mass of component $i$ per unit (current and referential respectively) volume occupied by component $i$;
    \item the volume fractions $\phi_i$ and $\phi_{0i}$ of phase $i$ in the current and reference respectively.
\end{enumerate}
We highlight that we will use a Lagrangian description for all of these quantities. That is, while $P_i, \rho_i, \phi_i$ all describe quantities in the current configuration, they will be written as functions of $\bfx_0$.

It follows from these definitions that:
\begin{subequations}
\label{eqn:true-density}
\begin{align}
    P_{0i} & = \phi_{0i} \rho_{0i}
    \\
    P_{i} & = \phi_{i} \rho_{i}
\end{align}
\end{subequations}
A central assumption is that the entire RVE is mapped by $\bfF$ as an affine transformation.
Further, the mass of each phase contained in an RVE is the same in the reference and current configurations\footnote{We note, however, that the mass of the fluid phases in both the reference and current configuration can change in time due to mass transport.}.
Together, these imply the following relations between quantities in the reference and current configurations:
\begin{subequations}
\label{eqn:ref-current-mapping}
\begin{align}
    & \text{Differential volume: } \dm\Omega_0 \mapsto \dm\Omega = J \dm\Omega_0
    \\
    & \text{Density: } P_{0i}  \mapsto P_{i} = J^{-1} P_{0i}
    \\
    & \text{True density: } \rho_{0i} \mapsto \rho_i = \frac{\phi_{0i}}{\phi_{i}} J^{-1} \rho_{0i}
\end{align}
\end{subequations}

The relation between the referential and current versions of other quantities such as the volume fraction will depend on assumptions such as incompressibility, and will be discussed in later sections.

We notice that $\phi_{0i}$ is a fixed referential quantity that depends only on the geometry of the reference RVE, while $\phi_i$  depends on the deformation.
Therefore, the latter is a suitable quantity to appear as an argument of the total free energy below.
Further, $P_{0i}$ is \textit{not} a fixed referential quantity, because the mass of the fluid phases in the RVE can change due to transport.
Therefore, $P_{0i}$ is a suitable quantity to appear as an argument of the total free energy.


\subsection{Energetics} \label{sec:Energetics}

The total free energy of the porous medium consisting of a solid skeleton and $N-1$ fluid phases is given by:
\begin{equation}
\label{eqn:energy}
    E\left[\bfx,\{P_{0i}\}, \{\phi_{i}\}\right] 
    = 
    \int_{\Omega_0} \left(W_0\left(\nabla_0\bfx,\{P_{0i}\}, \{\phi_{i}\}\right) -  \psi_0(\bfx, \{P_{0i} \}) \right) \dm\Omega_0
    -
    \int_{\partial A} \bft_0^* \cdot \bfx \dm S
    -
    \sum_{i=2}^N \int_{\partial B_i} \eta^*_{0i} P_{0i} \dm S
\end{equation}
where $W_0$ is the (Helmholtz) energy per unit referential volume, and $\psi_0$ is the force potential per unit referential volume. 
The boundary conditions correspond to specified mechanical traction $\bft^*$ on the solid skeleton over some part of the boundary $\partial A$, and specified referential chemical potential $\eta^*_{0i}$ over some part of the boundary $\partial B_i$ for each fluid phase.

We will assume that the energy density $W_0$ is additively composed of the energy densities of the components: 
\begin{equation}
\label{eqn:total-energy}
    W_0 = \phi_{0s} W_{0s}(\bfF) + \sum_{i=2}^N \phi_{0i} W_{0i}(J, P_{0i} , \phi_{i})
\end{equation}
where $W_{0i}$ is the energy density of phase $i$ per unit referential volume.

For the solid phase, it is standard to assume that $W_{0_1}$ is only a function of $\bfF$, and that $P_{0_1}$ and $\phi_1$ will not appear.

For the fluid phases, we are typically given a free energy density (per unit current volume) that is a function only of the true density of the corresponding phase.
Denoting this given free energy density as $W_i(\rho_i)$, we write the energy density $W_{0i}$ in terms of the arguments of $E$:
\begin{equation}
\label{eqn:fluid-energy}
    W_{0i} = \frac{J \phi_i}{\phi_{0i}} W_i \left(J^{-1} P_{0i} / \phi_i\right)
\end{equation}
where we have used that $\rho_i = J^{-1} P_{0i} / \phi_i$ from \eqref{eqn:true-density} and \eqref{eqn:ref-current-mapping}.

\subsection{Balance of Mass} \label{Balance of mass}

Following ~\cite{gajo2010general, coussy2004poromechanics}, we define the (referential) chemical potential $\eta_0$ as the variation of $E$ with respect to the referential density for the fluid phases in the bulk:
\begin{equation} \label{eq:chemical potential}
    \eta_{0i} 
    = 
    - \delta_{P_{0i}} E = - \parderiv{W_0}{P_{0i}} + \parderiv{\psi_0}{P_{0i}}
\end{equation}
where we use the notation $\delta_{P_{0i}} E$ to denote the variational (functional) derivative of $E$ with respect to $P_{0i}$.
The boundary condition that appears from the variational principle corresponds to $\eta_{0i} = \eta^*_{0i}$ specified on $\partial B_i$; this is related, but not directly equivalent, to specifying the pressure in the general setting.

The velocity in the reference configuration -- relative to the solid skeleton -- of the $i$-th fluid phase is assumed to be proportional to the gradient of the chemical potential:
\begin{equation} \label{eq:ref_relative velocity}
    \bfv_{0i} 
    = 
    \bfK_i \nabla_0\eta_{0i} 
    = 
    \bfK_i  \left(- \nabla_0 \parderiv{W_0}{P_{0i}} + \nabla_0 \parderiv{\psi_0}{P_{0i}} \right)
\end{equation}
where $\bfK_i$ is the positive-definite second-order permeability tensor of the porous medium in the reference configuration.
Following ~\cite{borja1995mathematical}, we can push forward \eqref{eq:ref_relative velocity} as follows:
\begin{equation} \label{eq:relative velocity}
    \bfv_{i} 
    = 
    J^{-1} \bfF \bfv_{0i} 
    = 
    \bfk_i  \left(- \nabla \parderiv{W_0}{P_{0i}} + \nabla\parderiv{\psi_0}{P_{0i}} \right)
\end{equation}
where $\bfk_i$ is the permeability tensor of the porous medium in the current configuration, $\bfK_i = J\bfF^{-1}\bfk_i \bfF^{-T}$ ~\cite{li2004dynamics}, using $\nabla_0 = \bfF^T \nabla$. In an isotropic medium, the permeability tensor can be written:
\begin{equation}
    \bfk_i = \frac{ \tilde{k}_i}{g} \bfI 
\end{equation}
where $g$ is the gravitational acceleration, and $\tilde{k}_i$ is the hydraulic conductivity, which is defined as $\tilde{k}_i = \frac{\kappa}{\gamma_i}\rho_i g$.
$\kappa$ is the true permeability of the solid; $\gamma_i$ is the dynamic viscosity of the fluid; $\rho_i$ is the fluid density; and $\bfI$ is the second-order identity tensor.

Therefore, the fluid flux vector for the $i$-th phase is:
\begin{equation} \label{eq:fluid flux}
    \bfq_i 
    = 
    P_i \bfv_i
    =
    -\bfk_i \left( P_i \nabla \parderiv{W_0}{P_{0i}} - P_i \nabla\parderiv{\psi_0}{P_{0i}} \right)
\end{equation}
The conservation of mass for each fluid phase in the current configuration can be written as follows:
\begin{equation}
\label{eqn:transport}
    -\int_{\partial\Omega} \bfq_i \cdot\bfn \dm S = \deriv{\ }{t} \left(\int_{\Omega} P_i \dm V\right)
    \Rightarrow
    -J \divergence \bfq_i = \deriv{\ }{t} P_i J
\end{equation}
where the differential form follows from the integral form by arbitrariness of $\Omega$. 

Specializing to a uniform gravitational field, the referential potential can be written as $\psi_0 = \bfb_0\cdot\bfx$, where $\bfb_0$ is the body force per unit referential volume.
This can be written in terms of the densities and volume fractions of the individual phases:
\begin{equation}
\label{eqn:grav-body-force}
    \bfb_0 := P_0\bfg 
    = \left( \sum_{i=1}^N P_{0i} \right) \bfg
\end{equation}
where $P_0$ is the total mass per unit referential volume, and $\bfg$ is the force per unit mass due to gravity. 
Then, we can simplify the potential term in \eqref{eq:fluid flux} to:
\begin{equation} \label{eq:flux gravity term}
    \parderiv{\psi_0}{P_{0i}} = \bfg\cdot\bfx 
    \Rightarrow
    \nabla\parderiv{\psi_0}{P_{0i}} = \bfg
\end{equation}
Consequently, the fluid flux vector in \eqref{eq:fluid flux}) can be rewritten as follows:
\begin{equation} 
    \bfq_i 
    = 
    -\bfk_i \left( P_i \nabla \parderiv{W_0}{P_{0i}} - P_i \bfg \right)
\end{equation}

\subsection{Balance of Momentum} \label{Balance of momentum}

Setting to zero the variation of $E$ with respect to $\bfx$ yields the force equilibrium equation in the referential and current forms:
\begin{equation}\label{eq:equlibrium in ref}
    \divergence_0 \bfT + \parderiv{\psi_0}{\bfx}=\mathbold{\num{0}} \Leftrightarrow \divergence \bfsigma + J^{-1} \parderiv{\psi_0}{\bfx}=\bf0
\end{equation}
where $\bfT := \parderiv{W_0}{\bfF}$ is the first Piola stress tensor; $\bfsigma = \frac{1}{J} \bfT \bfF^T$ is the Cauchy stress tensor; the body force is defined through the derivative of $\parderiv{\psi_0}{\bfx}$; and the corresponding boundary conditions from the variation are $\bfT \hat\bfn_0 = \bft_0^*$.

If we specialize to a uniform gravitational field, the body force in the reference $\bfb_0$ is given by the expression in \eqref{eqn:grav-body-force}.
In the current configuration, we have 
\begin{equation} 
    \bfb := P\bfg = \left(  \sum^N_{i=1} P_{i} \right) \bfg
\end{equation}
where $P$ is the total mass per unit current volume.
We note that $\bfb_0$ and $\bfb$ are related by the expression $\bfb = J^{-1} \bfb_0$.

\subsection{Balance of Fluid Pressure} \label{Balance of pressure}

We first notice that the variables $\phi_i$ are not independent when we take variations; they must satisfy the constraint $\sum^N_{i=1} \phi_i = 1$.
We can alternatively write this as:
\begin{equation}
    \sum^N_{i=1} \phi_i = 1 \Rightarrow \phi_2 = 1 - \phi_s - \phi_3 - \cdots - \phi_i - \cdots - \phi_N
\end{equation}
We will eliminate $\phi_2$ in terms of the other $\phi_i$.
Further, we notice that $\phi_s$ will be prescribed; similarly, some of the $\phi_i$ could be prescribed, e.g., due to incompressibility.

Setting to zero the variation of $E$ with respect to $\phi_i$, and using the form from \eqref{eqn:total-energy}, gives:
\begin{equation}
    \parderiv{W_0}{\phi_i} = \phi_{02} \parderiv{W_{02}}{\phi_2}\parderiv{\phi_2}{\phi_i} + \phi_{0i} \parderiv{W_{0i}}{\phi_i} = 0 
    \Rightarrow 
    - \phi_{02} \parderiv{W_{02}}{\phi_2} + \phi_{0i} \parderiv{W_{0i}}{\phi_i} = 0
\end{equation}
Since this is true for every $i$, it follows that:
\begin{equation}\label{eqn:pressure_equality}
    \phi_{02} \parderiv{W_{02}}{\phi_2} = \phi_{03} \parderiv{W_{03}}{\phi_3} = \cdots = \phi_{0N} \parderiv{W_{0N}}{\phi_N}
\end{equation}

We next recall that the derivative of the Helmholtz free energy with respect to volume, keeping temperature and mass fixed, is the pressure $p$.
Rewriting this in terms of the density $\rho$ -- which is inversely proportional to the volume when the mass is fixed -- and in terms of the Helmholtz free energy density per unit volume, we have the following relation between the fluid pressure $p$ and the Helmholtz free energy density:
\begin{equation}
    p = W(\rho) - \rho \deriv{W}{\rho}
\end{equation}
This relation holds for each fluid phase individually.

Applying the equilibrium relations above to the form assumed in \eqref{eqn:fluid-energy}, we get that:
\begin{equation}\label{eqn:phi_derivative}
    \phi_{0i}\parderiv{W_{0i}}{\phi_i}
    =
    J \left( W_i\left(J^{-1} P_{0i} / \phi_i \right) - \frac{J^{-1} P_{0i}}{\phi_i} \deriv{W_i}{\rho_i} \right)
    =
    J p_i
\end{equation}
Since $J$ depends only on $\bfx_0$, it follows that the pressures in each phase $p_i$ at a given $\bfx_0$ must be equal.

We highlight that the balance of pressure provides a local condition that must be satisfied pointwise and does not require the solution of PDE, in contrast to the balances of mass and momentum.

\subsection{Compatibility with the Second Law of Thermodynamics}
\label{2nd Law}

To ensure that the proposed model is consistent with the second law of thermodynamics, we first notice that in an isothermal setting it reduces to the requirement that the dissipation is non-negative for any process \cite{gurtin2010mechanics}.

Following \cite{penrose1990thermodynamically,abeyaratne2006evolution,marshall2014atomistic,de2018disclinations,agrawal2017dependence,agrawal2015dynamic,agrawal2015dynamic-2}, we begin by taking the time derivative of the free energy from \eqref{eqn:energy} to get:
\begin{equation}
    \frac{\dm E}{\dm t} = \int_{\Omega_0} \sum_i \left( \parderiv{W_0}{P_{0i}}  - \parderiv{\psi_0}{P_{0i}}  \right) \frac{\dm P_{0i}}{\dm t} \dm \Omega_0
\end{equation}
where we have used that $\deriv{\phi_i}{t} = 0$ and $\deriv{\bfx}{t} = \bf 0$ per our quasistatic model.
We have assumed that the boundary terms vanish for simplicity, but they can easily be accounted for following, e.g.,  \cite{penrose1990thermodynamically,agrawal2015dynamic}.

From Section \ref{Balance of mass}, we can push back the balance of mass to the reference configuration to obtain:
\begin{equation}
    -J \divergence \bfq_i = \frac{\dm}{\dm t} P_i J \Rightarrow -\divergence_0 \bfq_{0i} = \frac{\dm}{\dm t} P_{0i}
\end{equation}
and the corresponding flux for each fluid phase is given by:
\begin{equation}
    \bfq_{0i} = - P_i \bfK_i  \nabla_0 \left(   \parderiv{W_0}{P_{0i}} -  \parderiv{\psi_0}{P_{0i}} \right)
\end{equation}

Using these expressions, and imposing the condition that $\bfK$ is a positive-definite tensor, we can rewrite $\frac{\dm E}{\dm t}$ as follows and use integration-by-parts:
\begin{equation}
\begin{split}
    \frac{\dm E}{\dm t} &= 
        \int_{\Omega_0} \sum_i \left( \parderiv{W_0}{P_{0i}} - \parderiv{\psi_0}{P_{0i}} \right)  \divergence_0 \left[P_i \bfK_i \nabla_0 \left( \parderiv{W_0}{P_{0i}}- \parderiv{\psi_0}{P_{0i}} \right) \right] \dm \Omega_0 \\
        &= - \int_{\Omega_0} \sum_i P_i \nabla_0 \left( \parderiv{W_0}{P_{0i}}- \parderiv{\psi_0}{P_{0i}} \right) \cdot \bfK_i \nabla_0  \left( \parderiv{W_0}{P_{0i}} - \parderiv{\psi_0}{P_{0i}} \right) \dm \Omega_0 \\
        & \leq 0
\end{split}
\end{equation}
which shows that $E$ is non-increasing in the proposed model, as required by the second law of thermodynamics.


\section{Models for Saturated and Unsaturated Systems} \label{Constitutive Models for Saturated Systems}

\begin{proof}[Single Compressible Fluid in a Saturated Medium]
    The energy density is:
    \begin{equation}
        W_0(\bfF,P_{0f},\phi_{f})=(1-\phi_{0f})W_{0s}(\bfF)+\phi_{0f} W_{0f}(J,P_{0f},\phi_{f})
    \end{equation}
    where we have used $N=2$ (solid $s$ and fluid $f$).
\end{proof}

\begin{proof}[Multiple Immiscible Compressible Fluids]
    The energy density is:
    \begin{equation}\label{eqn:multiple fluid phases energy}
        W_0(\bfF, \{P_{0i}\}^N_{i=2}, \{\phi_{i}\}^N_{i=2}) 
        = 
        \phi_{0s} W_{0s}(\bfF) + \sum_{i=2}^N \phi_{0i} W_{0i}(J, P_{0i} , \phi_{i})
    \end{equation}
    where the subscript $s$ denotes quantities corresponding to the solid.
\end{proof}

\begin{proof}[Unsaturated Systems]
    Consider an unsaturated medium with a solid phase, a incompressible fluid phase (e.g. water), and an compressible fluid phase (e.g. air).
    The state variables for the energy are $\bfF, P_{0I}, P_{0c}, \phi_I, \phi_c$,
    where the subscripts $c$ and $I$ denote quantities corresponding to the compressible fluid and incompressible fluid respectively.
    
    For the incompressible fluid phase, we require that the true density in the current configuration has a fixed value $\tilde{\rho}_I$:
    \begin{equation}
        \frac{P_{0I}}{J \phi_I} = \tilde{\rho}_I
    \end{equation}
    We notice that the energy is constant for an incompressible fluid: all states that satisfy incompressibility have the same energy, and all other states are inadmissible (i.e., have infinite energy).
    Since this energy is a constant, it will not appear once we take variations and consequently we can set it to $0$ for simplicity.
    
    Therefore, we can write the total energy density as:
    \begin{equation}
        W_0(\bfF,P_{0I},P_{0c}, \phi_I, \phi_c)
        =
        \phi_{0s} W_{0s}(\bfF) 
        +  
        \phi_{0c} W_{0c}(J,\phi_c,P_{0c})
    \end{equation}
    We highlight some issues.
    First, if we further assume that the solid skeleton is incompressible and affinely deformed by $\bfF$, this imposes the additional constraint $J=1$.
    Second, we notice that the incompressibility constraint is generally not satisfied in the reference configuration, i.e. $\rho_{0I} \neq \tilde{\rho}_I$ in general.
    Third, it is sometimes reasonable to neglect the energy of the compressible fluid (e.g., air at low pressure in rock); in that case, we need only consider the energy of the solid.
    However, we emphasize that the problem still contains the influence of the fluid through the mass transport equations \eqref{eqn:transport}.
\end{proof}


\section{Models for Referential and Current Volume Fractions}\label{Models for Referential and Current Volume Fractions}

\begin{proof}[Incompressible Solid with Compressible Pore Space]

    \begin{figure}[htb!]
    	\includegraphics[width=0.8\textwidth]{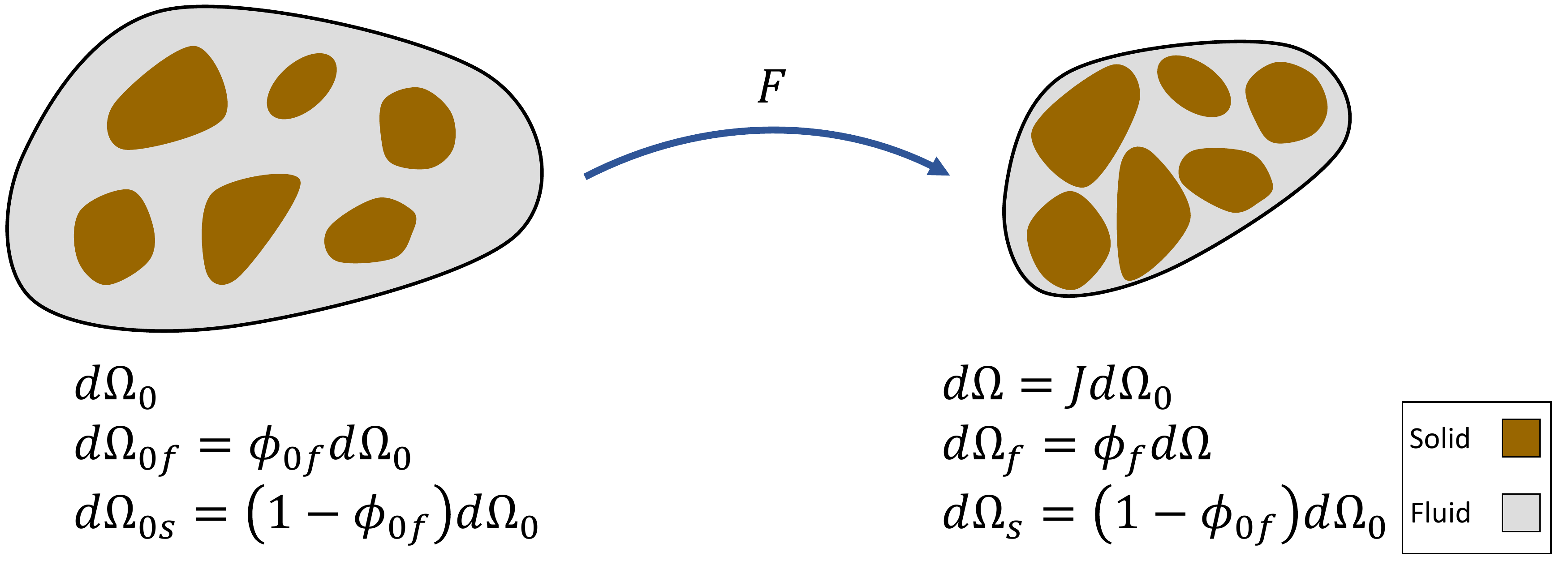}
    	\caption{Volume change for an incompressible solid phase with a compressible fluid.}
    	\label{fig:phi_1}
    \end{figure}

    Consider a porous medium with an incompressible solid skeleton and compressible fluid in the pore space (Figure \ref{fig:phi_1}).
    From incompressibility of the solid skeleton -- meaning that the true current density $\rho_s$ has a given fixed value -- and the fact that the mass of the solid skeleton is conserved because there is no mass transport of the solid, we have that the volume occupied by the solid phase in the reference and current is constant.
    Therefore, the current and referential volume fractions of the solid phase are related by $\phi_s = J^{-1} \phi_{0s}$.
    Using that the fluid volume fractions in the reference and current are given by $\phi_{0f} = 1 - \phi_{0s}$ and $\phi_{f} = 1 - \phi_s$, we have that $\phi_f = 1-J^{-1} (1-\phi_{0f})$.
\end{proof}

\begin{proof}[Affinely Deformed Solid with Compressible Fluid]

    \begin{figure}[htb!]
		\includegraphics[width=0.8\textwidth]{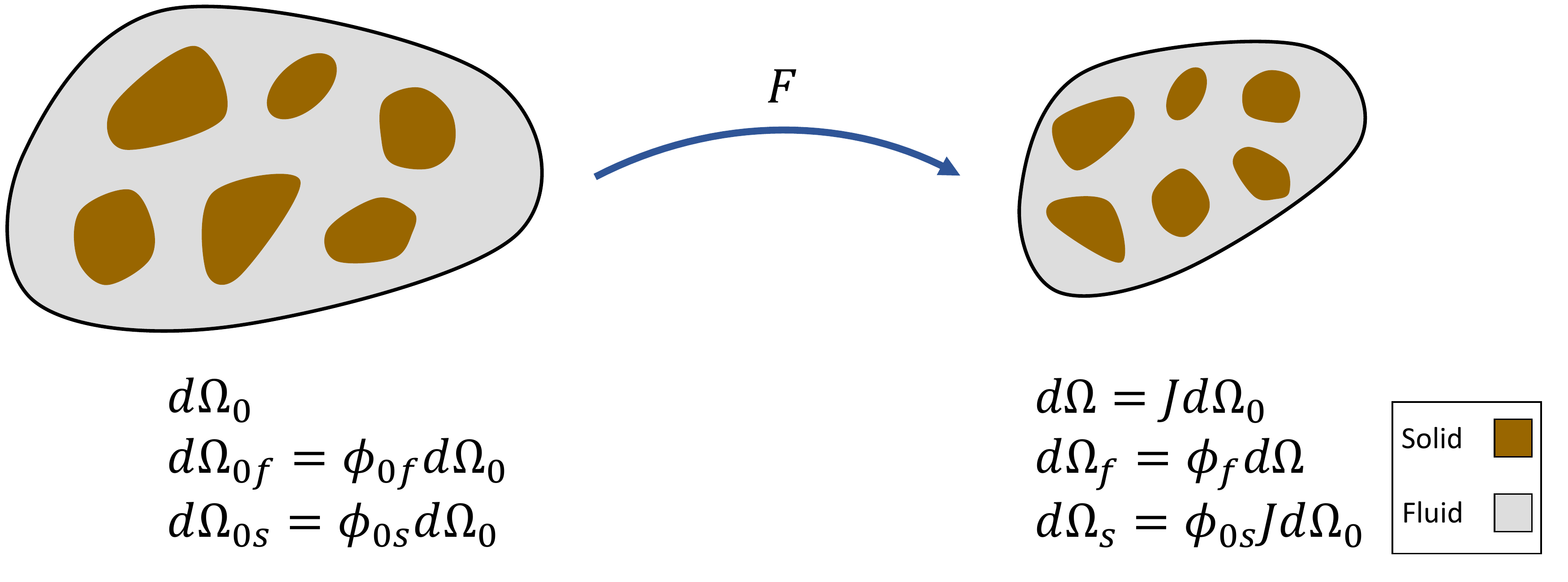}
        \caption{Volume change for a compressible solid phase with a compressible fluid.}
		\label{fig:phi_2}
    \end{figure}

    Consider a porous medium with a compressible solid phase and a compressible fluid phase (Figure \ref{fig:phi_2}).
    Let $\phi_{0f}$ and $\phi_{0s}$ be the referential volume fractions of fluid and solid phases, respectively. 
    The deformation of the solid phase is assumed to be affine with the macroscopic deformation gradient $\bfF$.
    This implies that $\phi_{0s} = \phi_s$ as follows.
    We have the following relation between current and reference volume elements for the entire RVE as well as for the solid skeleton:
    \begin{align}
        \dm \Omega = J \dm \Omega_0 \\ \nonumber
        \dm \Omega_s = J \dm \Omega_{0s}
    \end{align}
    Further, from the definition of the volume fraction of the solid skeleton, we have:
    \begin{align}
        \dm \Omega_s = \phi_s \dm \Omega = \phi_s J \dm \Omega_0 \\ \nonumber
        \dm \Omega_{0s} = \phi_{0s} \dm \Omega_0
    \end{align}
    therefore, by substituting last equation in second equation, we can derive:
    \begin{equation}
        \dm \Omega_{0} = J \phi_{0s} \dm \Omega_0
    \end{equation}
    We can therefore conclude that $\phi_{0s} = \phi_s$. 

    Further, since $\phi_f = 1 - \phi_s$ and $\phi_{0f} = 1 - \phi_{0s}$, it follows that $\phi_{0f} = \phi_f$.
\end{proof}

\begin{proof}[Compressible Solid with Compressible and Incompressible Fluids]

    \begin{figure}[htb!]
    	\includegraphics[width=0.8\textwidth]{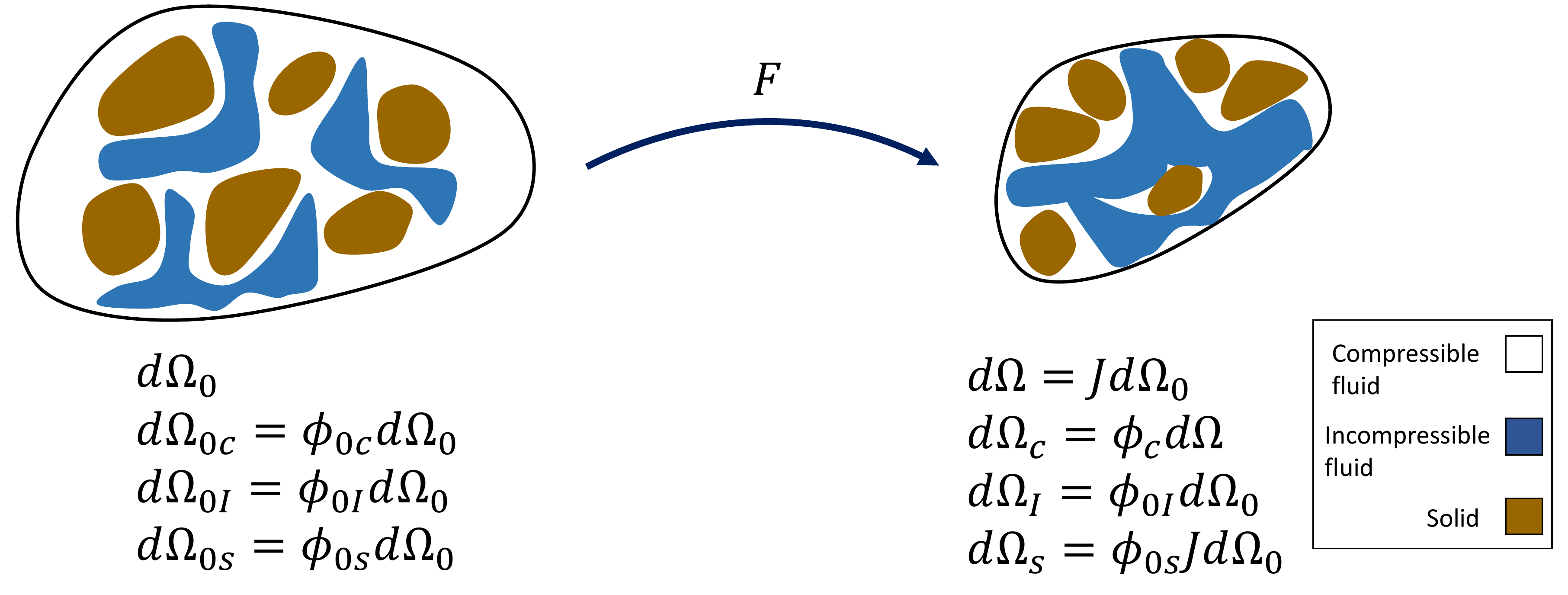}
        \caption{Volume change for a compressible solid phase with incompressible and  compressible fluid phases.}
    	\label{fig:phi_3}
    \end{figure}
    
    Consider a porous medium with compressible solid skeleton and air and incompressible water (Figure \ref{fig:phi_3}).
    Let $\phi_{0s}$, $\phi_{0I}$, and $\phi_{0c}$ be the reference volume fractions of compressible solid, incompressible fluid, and compressible fluid phases, respectively. 
    Assuming that the solid phase deforms affinely, we have that $\phi_s = \phi_{0s}$.
    For the incompressible fluid phase, we require that the true density in the current configuration have a fixed value $\tilde{\rho}_I$, giving $\phi_I = \frac{P_{0I}}{J \tilde{\rho}_I}$.
    Using that the volume fractions must sum to $1$, we have $\phi_c = 1 - \phi_{0s} - \frac{P_{0I}}{J \tilde{\rho}_I}$.
\end{proof}

\section{Model Verification with Terzaghi's and Mandel's Problems} 
\label{Model verification}

In this section, we verify the model by investigating its response in the model consolidation problems proposed by Terzaghi and Mandel. 
The general form of the strain energy density for a saturated porous medium with an incompressible solid phase and compressible fluid phase is:
\begin{equation}
    W_0 = (1-\phi_{0f})W_{0s}(\bfF) + \phi_{0f} W_{0f} (J, P_{0f}, \phi_f)
\end{equation}
Following section \ref{Balance of momentum}, the first Piola stress tensor can be derived as:
\begin{equation}
    \bfT = (1-\phi_{0f})\parderiv{W_{0s}}{\bfF} + \phi_{0f}\parderiv{W_{0f}}{J}J\bfF^{-T}
\end{equation}
Using \eqref{eqn:fluid-energy} and \eqref{eqn:phi_derivative}, we can write:
\begin{equation}
    \parderiv{W_{0f}}{J} = - p \frac{\phi_f}{\phi_{0f}}
\end{equation}
Thus, the Cauchy stress tensor can be written as follows:
\begin{equation}\label{eqn:Cauchy tensor linear}
    \bfsigma = (1-\phi_{0f})J^{-1}\bfF^{T}\parderiv{W_{0s}}{\bfF} - \phi_f p \bfI 
\end{equation}
Equation \eqref{eqn:Cauchy tensor linear} shows that the total Cauchy stress is obtained by a summation of the solid and fluid partial stresses $\bfsigma = \bfsigma_s + \bfsigma_f$, where $\bfsigma_f = -\phi_f p\bfI$ is an isotropic stress tensor. In addition, the general expression of the Cauchy effective stress tensor is $\bfsigma'=\bfsigma + p\bfI$, which can be written as:
\begin{equation}
    \bfsigma' = (1-\phi_f) \left( \bfF^T \parderiv{W_{0s}}{\bfF} +p \bfI\right)
\end{equation}
Assuming that the deformation of the solid phase is small and that the response is isotropic, we linearize it to get:
\begin{equation}
    \bfsigma = (1-\phi_{f})\left( 2\mu\bfepsilon + \lambda \trace(\bfepsilon)\bfI \right) - \phi_f p \bfI 
\end{equation}
where $\trace$ denotes the trace operator.
We then define the renormalized Lame constants $\tilde{\mu} = (1-\phi_{0f})\mu$ and $\tilde{\lambda} = (1-\phi_{0f})\lambda$. 
We define $\tilde{p} = \phi_f p$ and $p$ as the total pressure and the intrinsic fluid pressure, respectively. 

Following section \ref{Balance of mass}, the conservation of mass for the fluid phase can be written as:
\begin{equation} 
    \frac{\dm }{\dm t}P_f J = -J \divergence \bfq
\end{equation}
Following \cite{borja2005conservation}, the bulk modulus of the fluid phase can be written:
\begin{equation} \label{eq:bulk modulus}
    K_f = \rho_f \deriv{{p}}{\rho_f}
\end{equation}
We recall further that:
\begin{equation} \label{eq:m and rho relation}
    P_f = \phi_f\rho_f
\end{equation}
Using \eqref{eq:bulk modulus} and \eqref{eq:m and rho relation}, we can rewrite the conservation of mass for the fluid and solid phases respectively as:
\begin{align} 
    \label{eq:6.10}
    & \frac{\dm J}{\dm t}\frac{\phi_f}{J} + \frac{\dm \phi_f}{\dm t} + \frac{\dm p}{\dm t}\frac{\phi_f}{K_f} = -\frac{1}{\rho_f}\divergence \bfq
    \\
    & \label{eq:6.11}
    \frac{\dm J}{\dm t}\frac{\phi_s}{J} + \frac{\dm \phi_s}{\dm t} = 0   
\end{align}
The second equation, for the solid phase, is derived in the same way as the equation \eqref{eq:6.10} for the fluid phase, but also using that there is no transport of the solid phase and that the solid phase is incompressible (the bulk modulus $K_s\to\infty$).

Adding \eqref{eq:6.10} and \eqref{eq:6.11} yields the balance of mass of the porous system as follows:
\begin{equation}\label{eq:total mass conservation}
    \frac{\dm J}{\dm t} \frac{1}{J} + \frac{\dm p}{\dm t}\frac{\phi_f}{K_f} = -\frac{1}{\rho_f}\divergence \bfq
\end{equation}
Ignoring the effect of gravity and using \eqref{eqn:ref-current-mapping} and \eqref{eq:chemical potential}, it follows that:
\begin{align}
    & \eta_0 = \parderiv{W_f(\rho)}{\rho}
    \\\label{eqn:5.15}
    & \parderiv{p}{\eta_0} = -P_f
\end{align}
Using these equations and considering small deformations, \eqref{eq:total mass conservation} can be written as:
\begin{equation}
    \frac{\dm \trace(\bfepsilon)}{\dm t} + \frac{\dm p}{\dm t}\frac{\phi_f}{K_f} = \frac{\kappa}{\gamma}\nabla^2 {p}
\end{equation}

\subsection{Terzaghi's one-dimensional consolidation} \label{Terzaghi's one dimensional consolidation}

\begin{figure}[htb!]
	\begin{center}
		\includegraphics[width=0.45\textwidth]{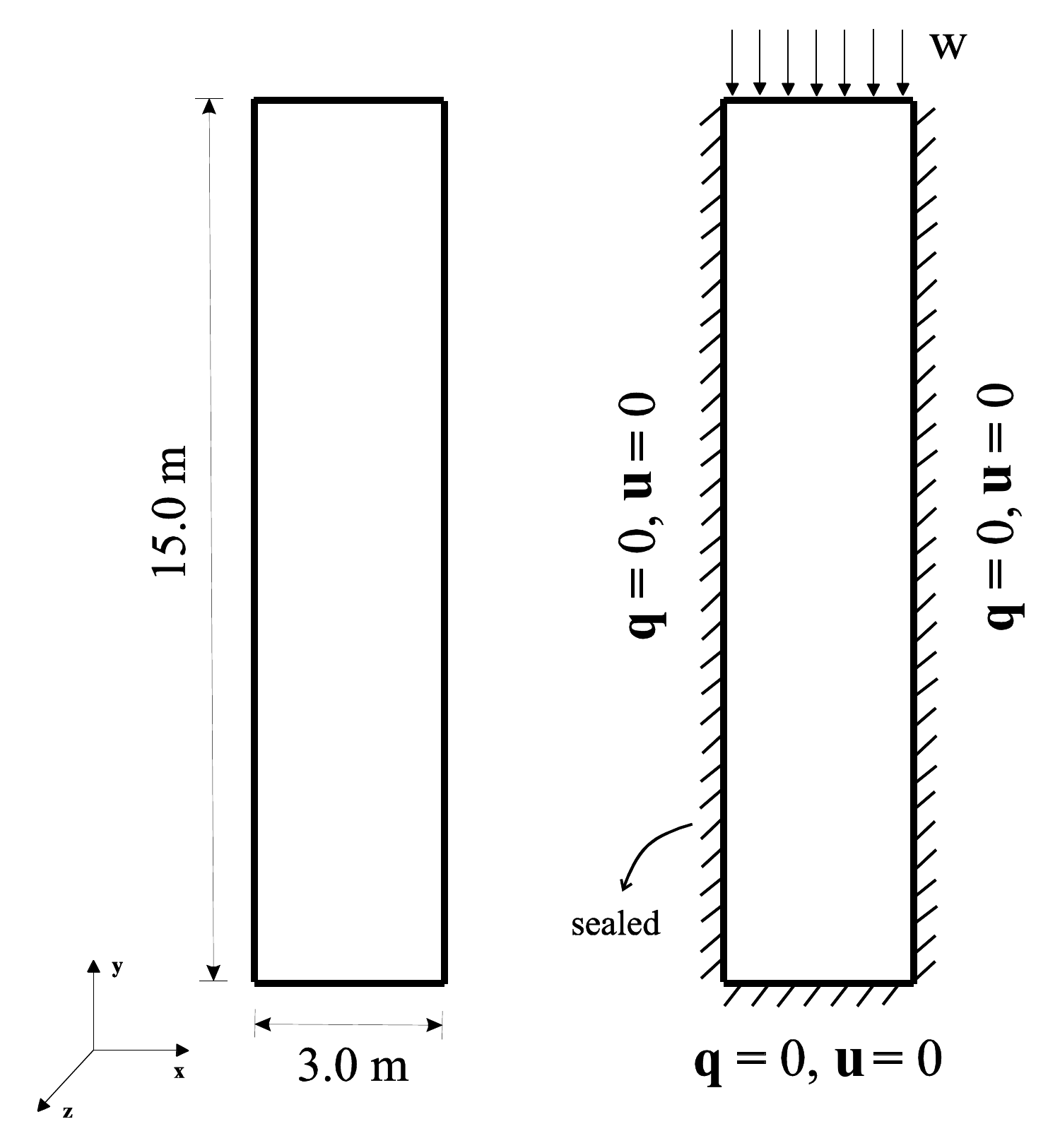}
		\caption{Geometry of Terzaghi's problem.}
		\label{fig:2}
	\end{center}
\end{figure}

In this section, we investigate the one-dimensional Terzaghi's problem. 
The geometry of the sample is shown in Fig. \ref{fig:2}, and the numerical values of the model parameters are listed in Table \ref{tab:1}. 
We consider a fully saturated sample with porosity $\phi_f = \num{0.375}$ that is loaded instantaneously at $t=\num{0} $ by a vertical load $w = \SI{10}{\kilo\pascal}$ at the top boundary, and then remains constant in time.
This sets up a mechanical stress in the sample.

The fluid flux is allowed to drain out of the top surface and is zero on the bottom surface.
The pressure boundary conditions of are:
\begin{equation}
    \text{ Top BC: } \tilde p(z=0) = 0; \quad \text{ Bottom BC: } \left.\parderiv{\tilde p}{z}\right|_{z=h} = 0
\end{equation}

The closed-form solution to this problem can be written as \cite{coussy2004poromechanics}:
\begin{align}
    \Bar{p} (\Bar{z},\Bar{t}) = \sum_{n=0}^{\infty} \frac{4}{\pi (2n+1)} \sin \left( \frac{(2n+1)\pi}{2}\Bar{z} \right) \exp \left( -\frac{(2n+1)^2 \pi^2}{4} \Bar{t} \right)
    \\
    \tilde{p} = \frac{w K_f }{K_f + \phi_f (\tilde{\lambda} + 2\tilde{\mu})}\Bar{p}; \quad z = h\Bar{z}; \quad t = \frac{h^2\kappa}{\gamma}\left( \frac{ K_f (\tilde{\lambda} + 2\tilde{\mu})}{K_f + \phi_f (\tilde{\lambda} + 2\tilde{\mu})} \right) \Bar{t}
\end{align}

\begin{table*}[h]
    \centering
    \begin{tabular}{l l}
        \hline
        Property & Value\\
        \hline
        Solid phase Lame constant, $\tilde{\lambda}$ \qquad \qquad & \SI{40}{\mega\pascal}\\
        Solid phase Lame constant, $\tilde{\mu}$ & \SI{40}{\mega\pascal}\\
        Hydraulic conductivity, $\tilde{k}$ & \SI{10e-8}{m\per s} \\
        Fluid bulk modulus, $K_f$ &  \SI{2270}{\mega\pascal}\\
        \hline
    \end{tabular}
    \caption{Properties of solid and fluid phases.}\label{tab:1}
\end{table*}

\begin{figure}[htb!]
	\begin{center}
		\includegraphics[width=0.8\textwidth]{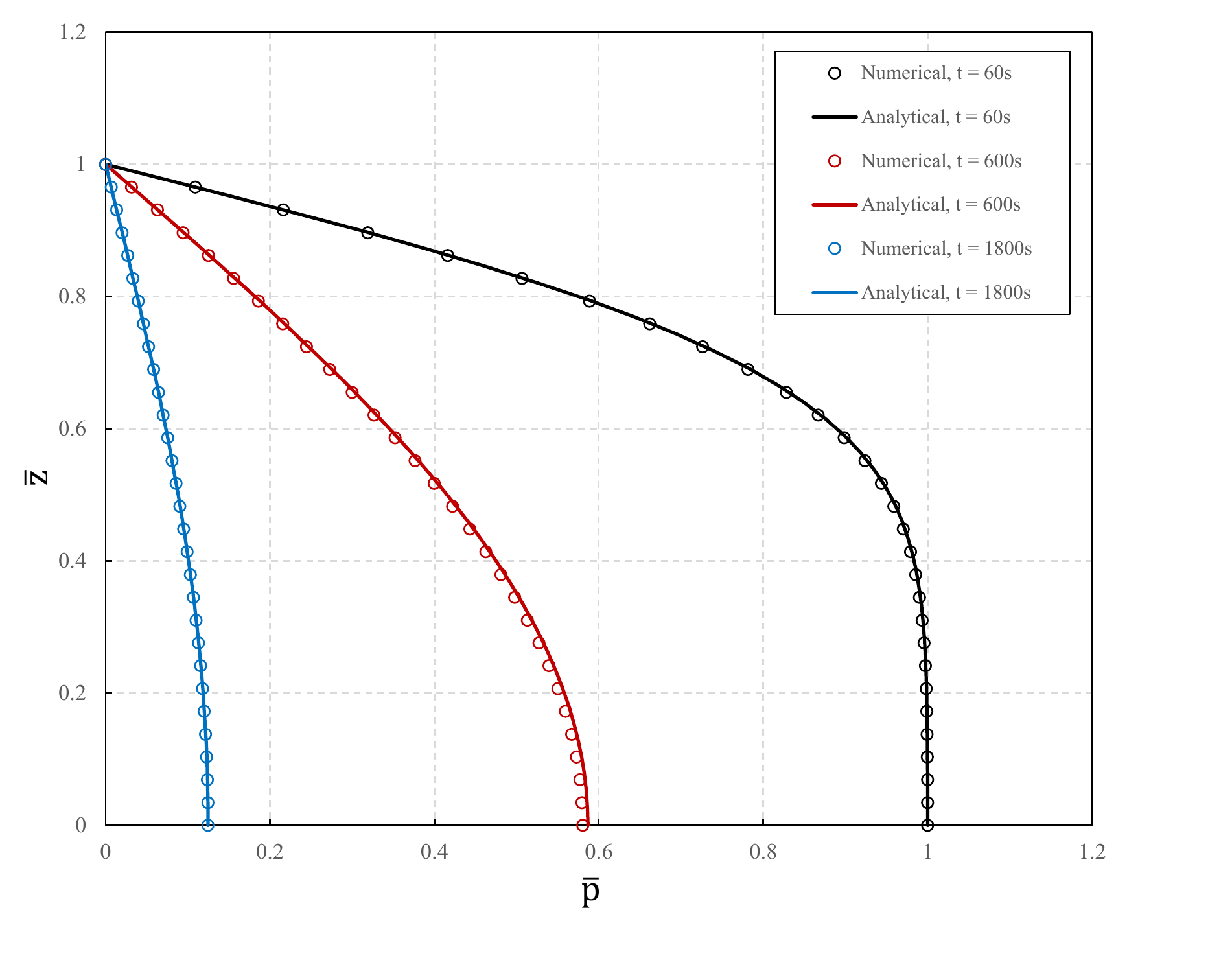}
		\caption{Variation of pressure with height for Terzaghi's problem. }
		\label{fig:4}
	\end{center}
\end{figure}

Fig. \ref{fig:4} shows the comparison of the closed-form and numerical results at different times.
Due to the compressibility of the pore space, the fluid drains gradually from the top boundary of specimen; therefore the pore pressure value decreases and the effective stress increases as time increases.
As time increases, the pore pressure tends to $0$ and the solid partial stress tends to the external load.

\subsection{Mandel's two-dimensional consolidation} \label{Mandel's two dimensional consolidation}

\begin{figure}[htb!]
	\begin{center}
		\includegraphics[width=0.8\textwidth]{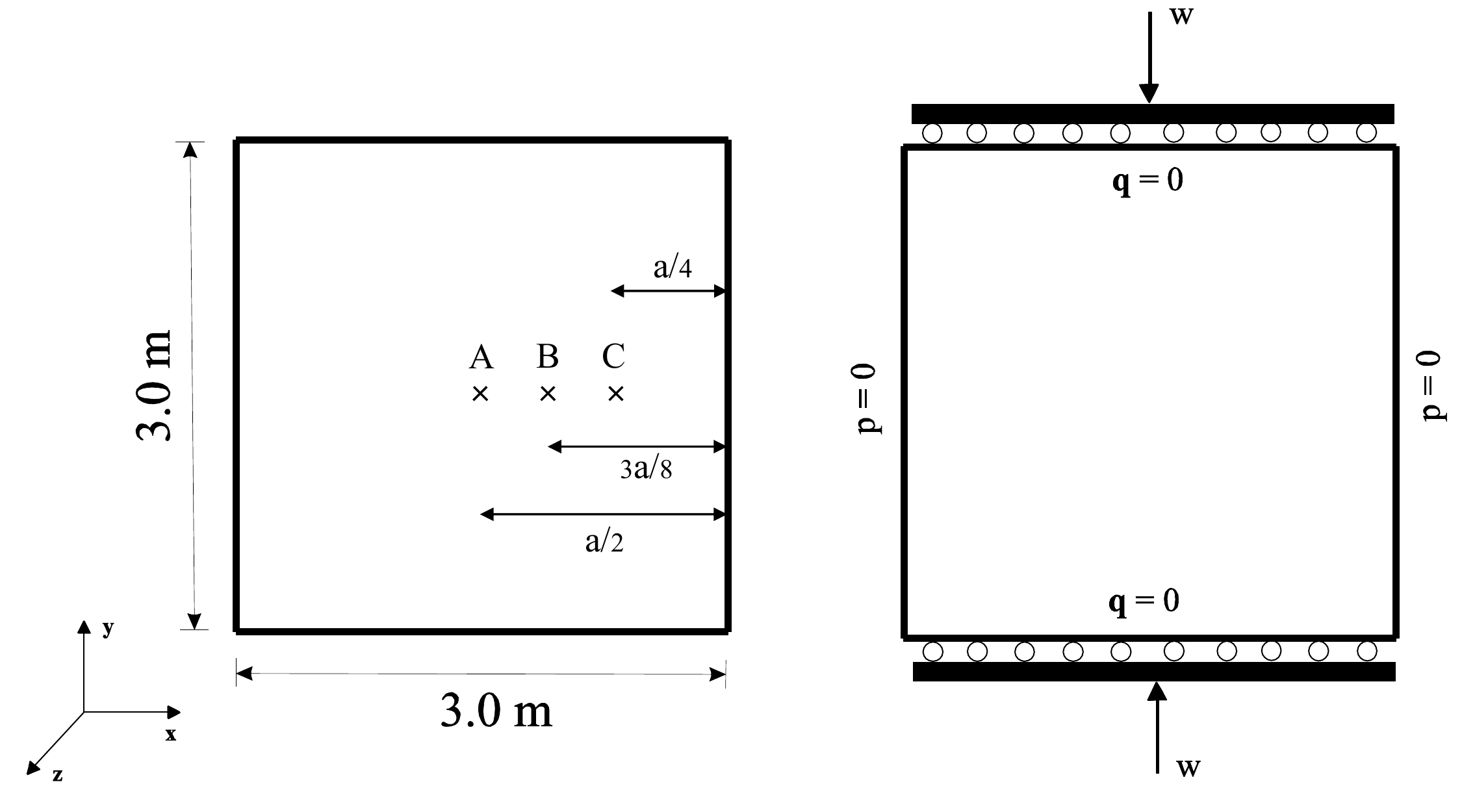}
		\caption{Geometry and the boundary conditions for Mandel's problem.}
		\label{fig:6}
	\end{center}
\end{figure}

In this section, we discuss a two-dimensional consolidation of a porous sample. 
The geometry of the specimen is shown in Fig. \ref{fig:6}, and the material properties and soil porosity are assumed the same as the previous example. 
We consider a specimen with infinite length in the $z$ direction which is between two horizontal plates (the top and bottom boundaries).

A constant load $w = \SI{10}{\kilo\pascal}$ is applied in the vertical direction at $t = \num{0}$, and remains constant thereafter. 
The left and right boundaries are traction-free.

For the fluid flow, we assume that there is no flux at the top and bottom boundaries, and the pressure is $0$ on the left and right boundaries (located at $x=\pm a/2$).

The closed-form solution for the pressure distribution is \cite{coussy2004poromechanics}:
\begin{align}
    \Bar{p}(\Bar{x},\Bar{t}) = 2 \sum_{n=1}^{\infty} \frac{\cos (\alpha_n \Bar{x})-\cos \alpha_n}{\alpha_n - \sin \alpha_n \cos \alpha_n}\sin\alpha_n \exp (-\alpha_n^2 \Bar{t})
    \\
    \tilde{p} = \frac{1}{3}B(1+\nu_u )w\Bar{p} ; \quad x = a\Bar{x} ;\quad t = \frac{a^2\kappa}{\gamma}\left( \frac{ K_f (\tilde{\lambda} + 2\tilde{\mu})}{K_f + \phi_f (\tilde{\lambda} + 2\tilde{\mu})} \right) \Bar{t}
\end{align}
where $\alpha_n$ is the solution of the following equation:
\begin{equation}
    \frac{\tan \alpha_n}{\alpha_n}=\frac{1-\nu}{\nu_u - \nu}; \quad \text{ where }
    \nu = \frac{\tilde{\lambda}}{2(\tilde{\lambda} + \tilde{\mu})} ; \quad \nu_u = \frac{3\nu + B(1-2\nu)}{3-B(1-2\nu)} ; \quad B= \frac{K_f}{\phi_f K_u} ; \quad K_u = \lambda+\frac{2\tilde{\mu}}{3}+\frac{K_f}{\phi_f}
\end{equation}

\begin{figure}[htb!]
	\begin{center}
		\includegraphics[width=0.65\textwidth]{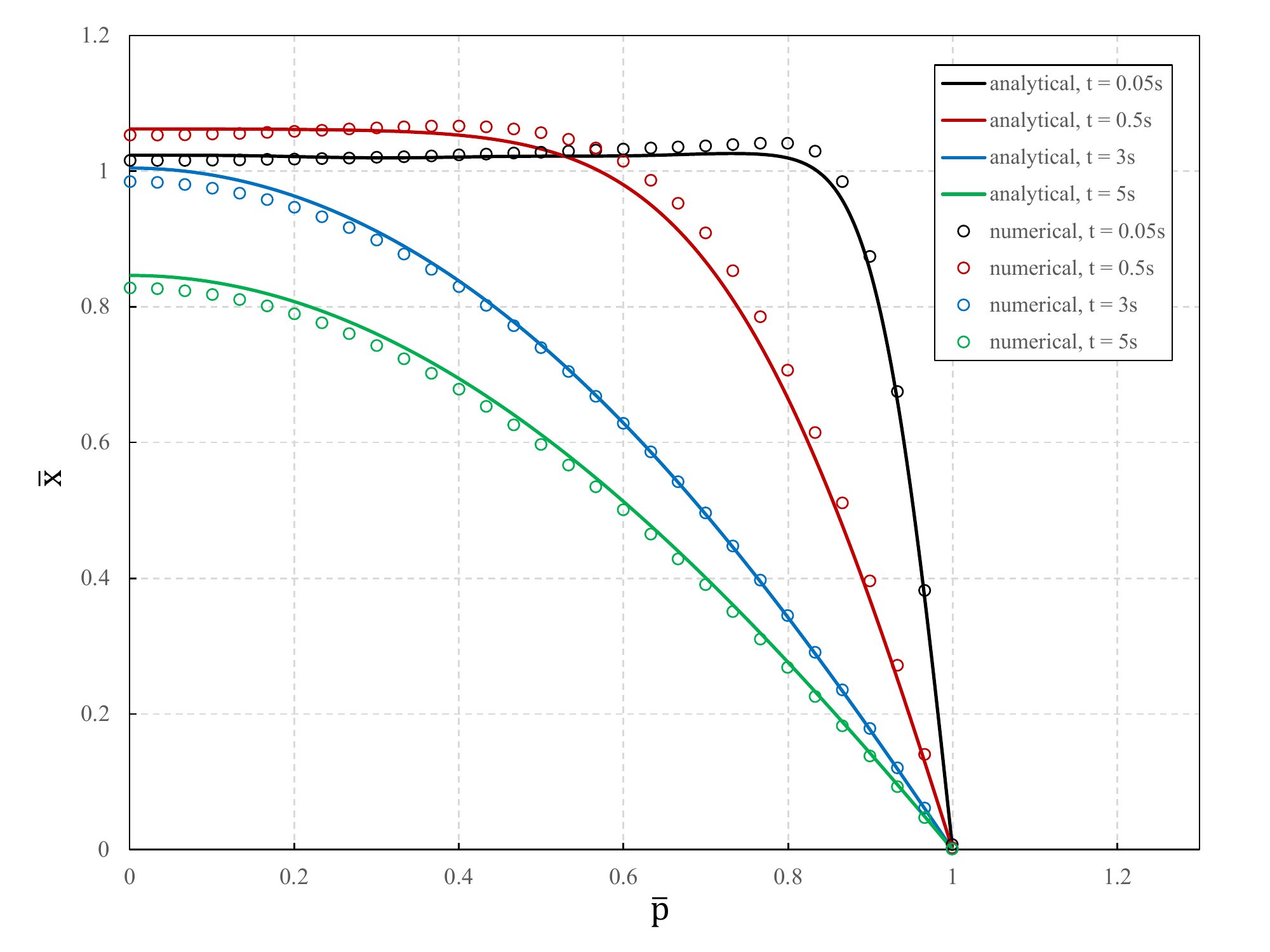}
		\caption{Variation of pore pressure in the horizontal direction at different times.}
		\label{fig:7}
	\end{center}
\end{figure}

\begin{figure}[htb!]
	\begin{center}
		\includegraphics[width=0.55\textwidth]{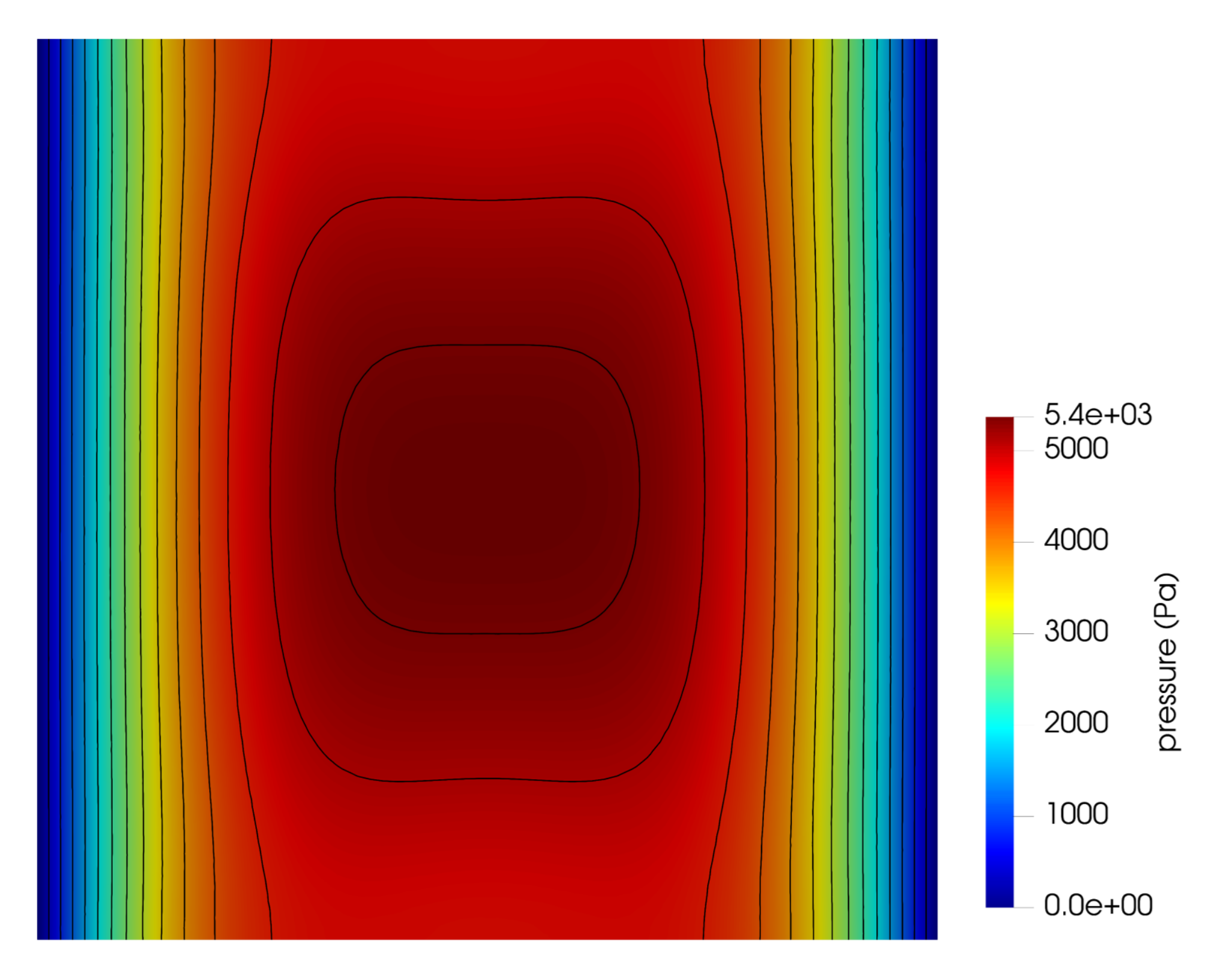}
		\caption{Pressure distribution for Mandel's problem at $t=\SI{1}{\second}$.}
		\label{fig:8}
	\end{center}
\end{figure}

Once the load is applied, the fluid pressure decreases at the vicinity of left and right boundaries due to drainage out of the sides.
Due to relatively slow drainage, the pressure decrease cannot be observed immediately in the whole domain, while the total load remains constant. 
Therefore, the pore pressure increases beyond the initial value at the center of specimen, which is the Mandel effect. 
The closed-form solution and our numerical results are shown in Fig. \ref{fig:7} for various times; we notice a signature of the Mandel effect at $t = \SI{0.5}{\second}$ wherein we have the maximum value of pore pressure in the central region of the specimen. 
Similarly, Fig. \ref{fig:8} shows the pore pressure distribution at $t=\SI{1.0}{\second}$.
We see clearly the increase of pore pressure beyond the initial pressure at the center, and it goes to $0$ on the right and left boundaries. 

\begin{figure}[htb!]
	\begin{center}
		\includegraphics[width=0.65\textwidth]{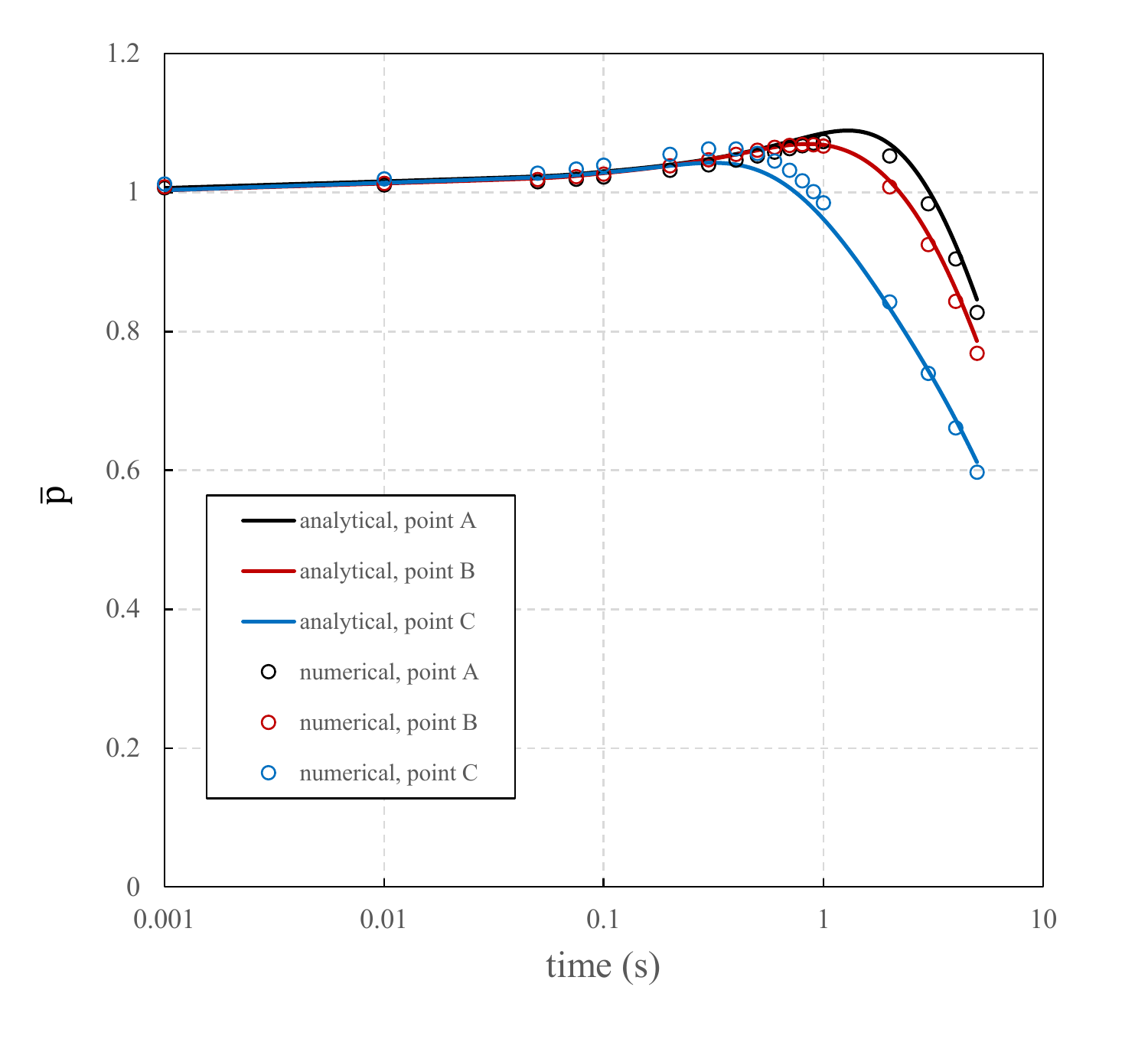}
		\caption{Pore pressure vs. time at different distances from the center of specimen.}
		\label{fig:9}
	\end{center}
\end{figure}

The variation of non-dimensional pore pressure ($\Bar{p}$) at different locations (from Fig. \ref{fig:6}) vs. logarithmic time is plotted in Fig. \ref{fig:9}. 
The log scale makes more prominent the pore pressure increase above its initial value shortly after applying the vertical load due to the Mandel effect.
Also, the comparison of these three graphs at points $A$, $B$, and $C$, clearly shows that the specimen experiences the highest amount of early times pressure increase at the central points (point $A$), and the Mandel effect vanishes at the vicinity of left and right boundaries.

\section{Van der Waals compressible fluid} \label{single compressible fluid}

In this section, we examine a porous medium containing a single fluid phase that is modeled as a Van der Waals (VdW) gas, motivated by experimental approaches such as \cite{krevor2012relative}.
A key feature of the VdW model is that it allows for a change from gas to liquid, and provides a simple model for gases such as carbon dioxide (CO2).
An important challenge in simulating gas injection in porous media is to predict the gas-liquid phase transformation of the injected gas under high pressure in the vicinity of the injection point.
The VdW model for the fluid phase enables us to investigate this complex behavior.

The Helmholtz free energy density per unit current volume of a VdW gas is given by the expression:
\begin{equation}
    W_f(\rho) = c \rho RT\left( 1-\log(cRT)\right)-\rho RT\log\left( \rho-b\right) -a\rho^2
\end{equation}
where $c$, $a$, $b$, $R$ and $S_0$ are constants, and $T$ is the temperature \cite{muller2009fundamentals}.

Following \eqref{eqn:fluid-energy}, we rewrite this as:
\begin{equation}
    W_{0f}(J,P_0,\phi)
    =
    \frac{J \phi}{\phi_0} W\left(J^{-1} P_0 / \phi\right)
    =
    c \frac{P_0}{\phi} RT\left( 1-\log(cRT)\right)-\frac{P_0}{\phi} RT\log\left( \frac{\phi J}{P_0}-b\right) -a\frac{P_0^2}{\phi^2 J}
\end{equation}
where we have used that $\phi=\phi_0$ from Section \ref{Models for Referential and Current Volume Fractions} (affinely deformed solid with compressible fluid).

We use a Neo-Hookean model for the solid phase:
\begin{equation} \label{eq:Neo-Hookean}
    W_{0s}(\bfF)
    =
    \frac{\mu}{2} (\trace(\bfF^T \bfF)-2) 
    - \mu \log J
    + \frac{\lambda}{2}(\log J)^2 
\end{equation}
Using \eqref{eq:chemical potential} and \eqref{eq:fluid flux}, we find the chemical potential and fluid flux:
\begin{align}
    \eta_0 & = \left( cRT\left( 1-\log(cRT)\right) -RT\log\left(\frac{1}{\rho} -b\right)+\frac{RT}{1-b \rho} -2 a \rho \right)
    \\
    \bfq & = -\phi\bfk \left( \frac{RT\nabla \rho}{1-b\rho} + \frac{bRT\rho\nabla \rho}{(1-b\rho)^2} -2 a\rho\nabla \rho \right)
\end{align}
where $\rho = \frac{P}{\phi}$ is the true current density of the fluid.

The first Piola stress can be obtained from the energy density:
\begin{equation}
    \bfT = (1-\phi)\left(\mu\bfF - \mu\bfF^{-T} +\lambda\log{J}\bfF^{-T} \right) - P_0\left( \frac{RT \bfF^{-T}}{1-b\rho} - a\rho\bfF^{-T} \right)
\end{equation}

\begin{figure}[htb!]
	\begin{center}
		\includegraphics[width=0.75\textwidth]{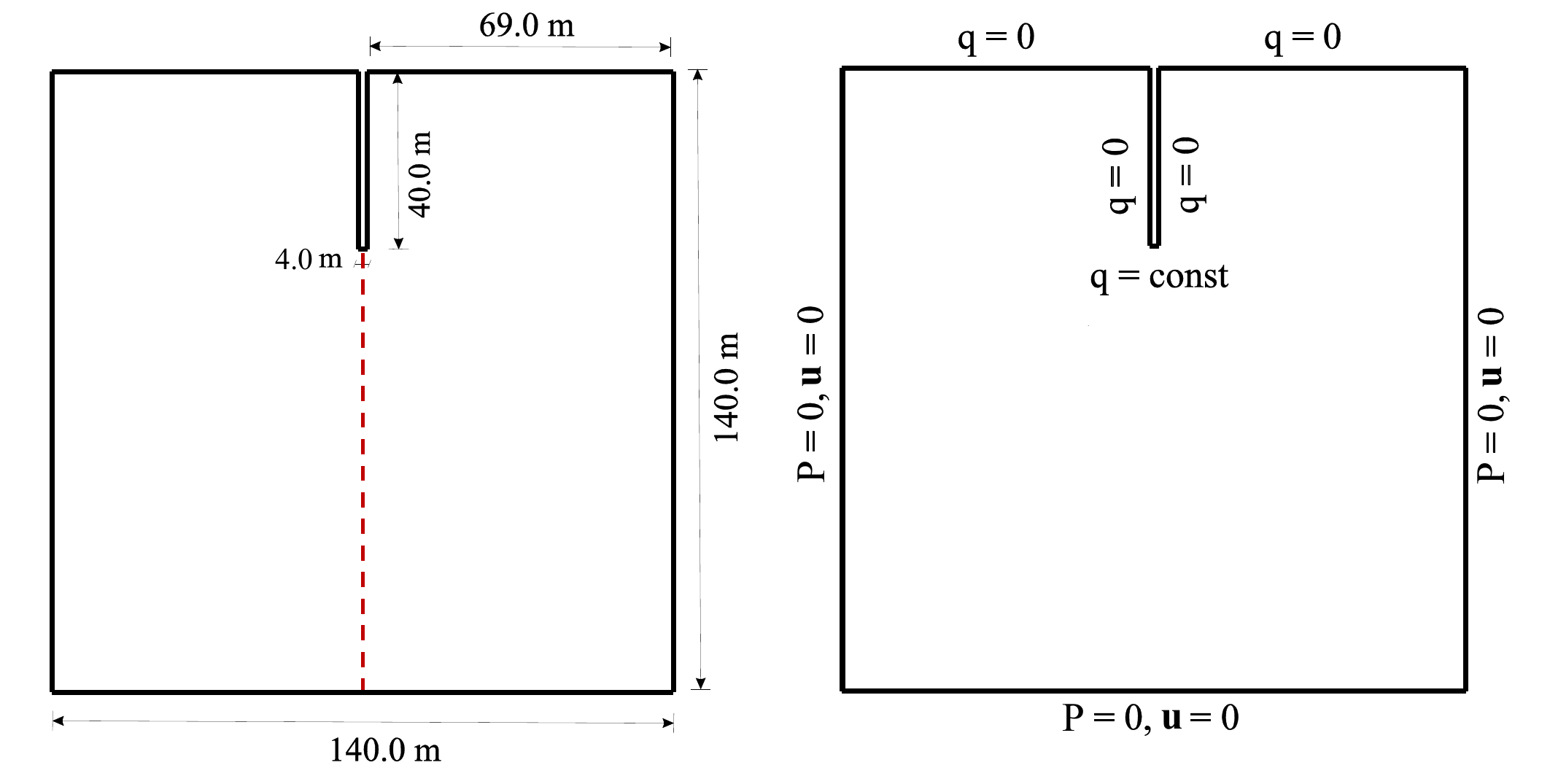}
		\caption{Geometry and boundary conditions of the specimen.}
		\label{fig:10}
	\end{center}
\end{figure}

The geometry and boundary conditions are shown in Fig. \ref{fig:10}, and the material properties of soil and CO2 gas are listed in table \ref{tab:2}.
This provides a simplified model of CO2 injection through a well.
The numerical calculation uses soil porosity of $\phi_{0f} = \num{0.2}$ and a constant flux $\bfq = \SI{4.4}{\kilo\gram \per \meter^2.\second}$ $\left(\SI{2.225e3}{\liter \per \meter^2.\second} \right)$ is assumed at the injection point. 

\begin{table*}[h]
    \centering
    \begin{tabular}{l l}
        \hline
        Property & Value\\
        \hline
        Solid phase Lame constant, $\lambda$ \qquad \qquad & \SI{57.7}{\mega\pascal}\\
        Solid phase Lame constant, $\mu$ & \SI{38.46}{\mega\pascal}\\
        Hydraulic conductivity, $\tilde{k}$ & \SI{0.1}{m\per s}\\
        Gas constant, $R$ & \num{8.32} \si{{\cubic\metre.\pascal}\per {\kelvin.\mole}} \\
        CO$_2$ constant, $a$ & \num{.364} \si{\pascal.\metre^6\per\mole^2} \\
        CO$_2$ constant, $b$ & \num{42.67e-6} \si{\metre^3\per\mole}\\
        Critical temperature, $T_c$  & \SI{303.4}{\kelvin}\\
        \hline
    \end{tabular}
    \caption{Properties of the solid phase and VdW gas (corresponding to CO2, from \cite{bielinski2007numerical})}.\label{tab:2}
\end{table*}

\begin{figure}[htb!]
	\begin{center}
		\includegraphics[width = 0.9 \textwidth]{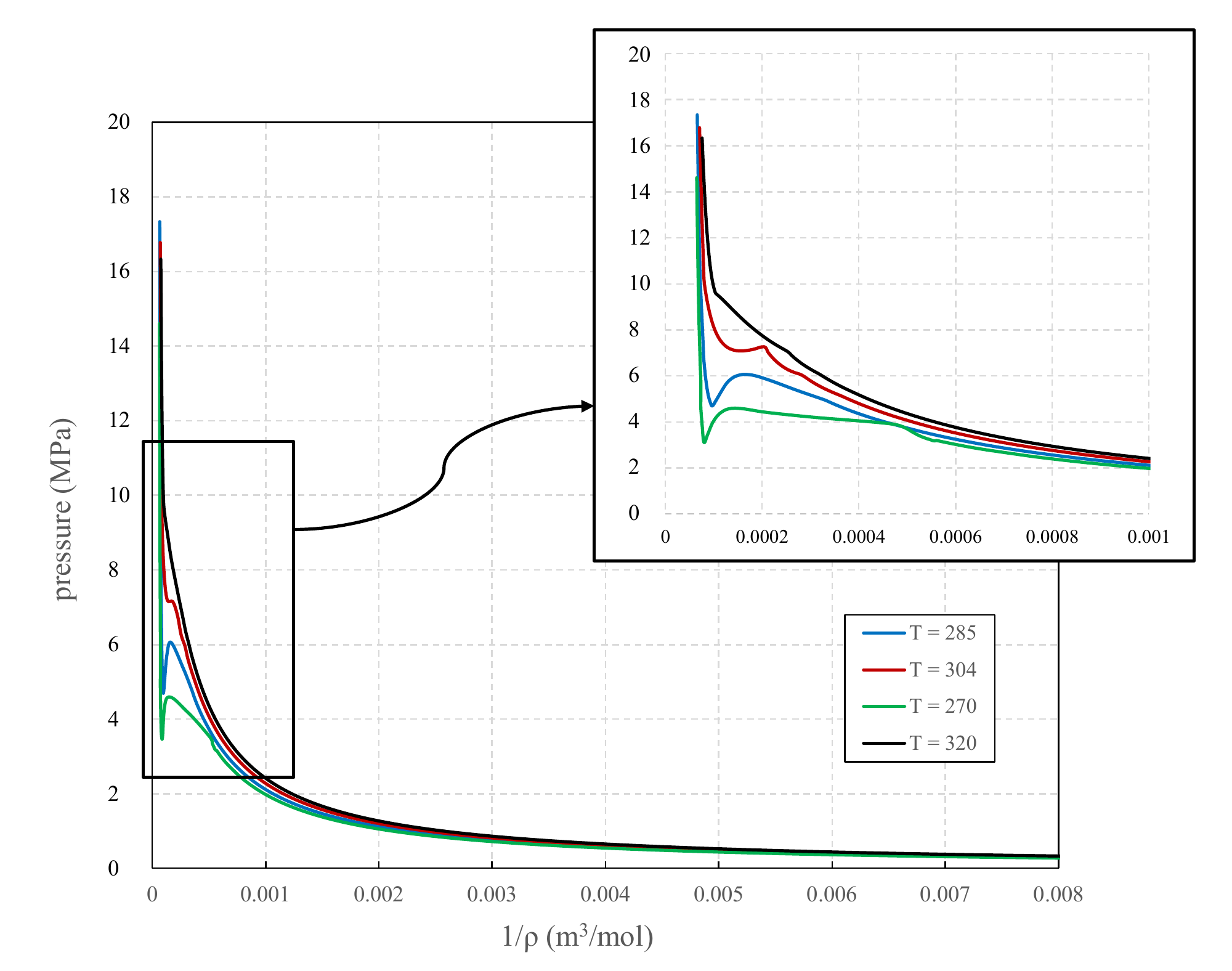}
		\caption{Pore pressure variation vs. inverse fluid density at different temperatures at $t = \SI{20}{\second}$.}
		\label{fig:12}
	\end{center}
\end{figure}

Fig. \ref{fig:12} illustrates the variation of pore pressure vs. inverse fluid density ($\rho^{-1}$), along the vertical direction below the center of injection (shown by the dashed line in Fig \ref{fig:10}). 
As shown in Fig. \ref{fig:12}, the simulation is repeated at different temperatures: below, above, and close to the critical temperature of CO2.
As expected, near the center of injection, shortly after start of fluid injection, we see a gas-liquid phase transformation due to the high pressures when the temperatures are below the critical temperatures ($T < T_c, ~ T= \SI{270}{\kelvin}$ and $T= \SI{285}{\kelvin}$).

\begin{figure}[htb!]
	\begin{center}
		\includegraphics[width=0.9\textwidth]{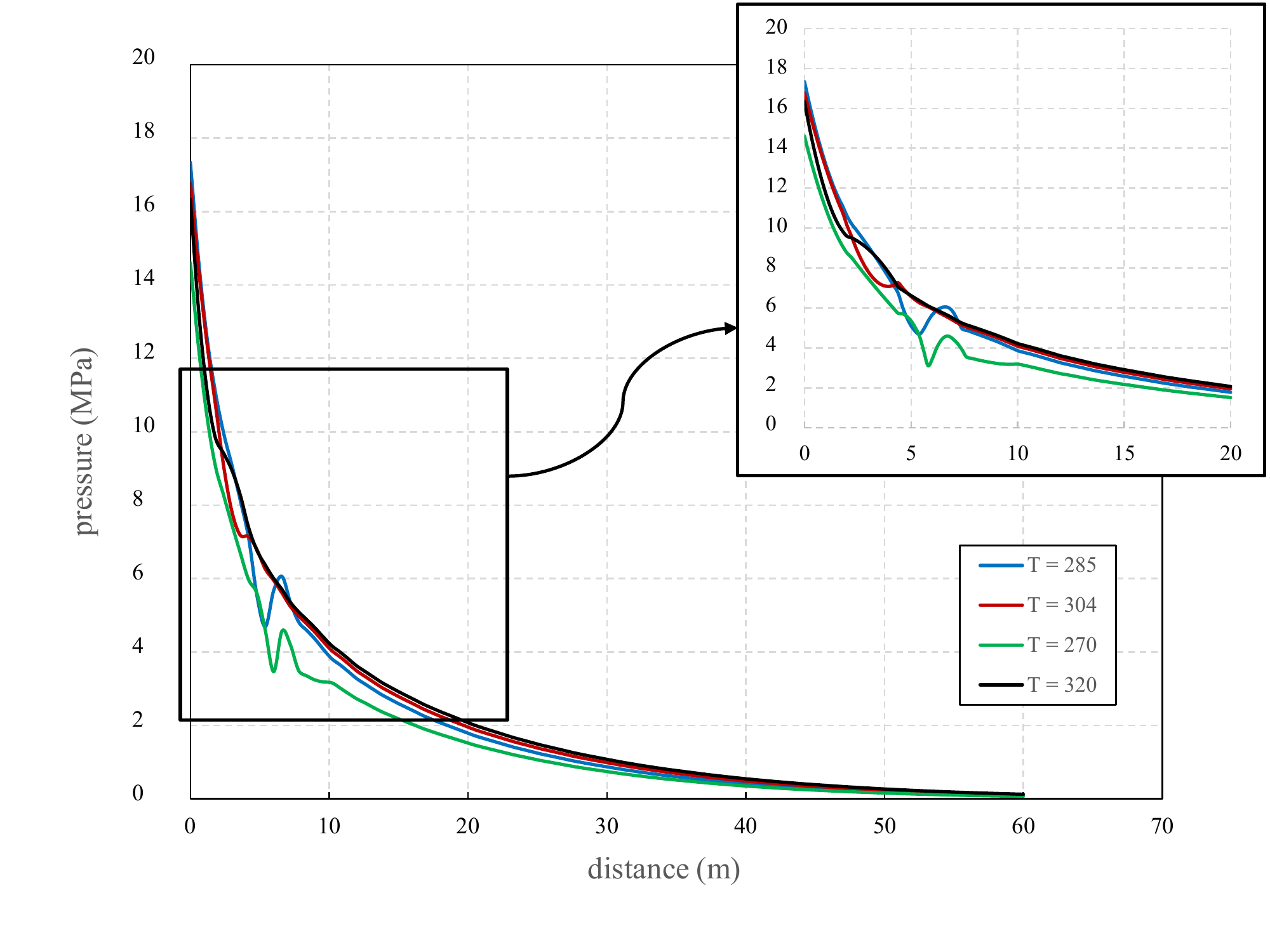}
		\caption{Pore pressure variation vs. distance from the injection point at $t=20s$ for different temperatures.}
		\label{fig:13}
	\end{center}
\end{figure}

Fig. \ref{fig:13} plots the variation of pore pressure with distance from the injection point along the vertical direction.
Near the center of injection, we have the highest value of pore pressure which decreases gradually as we move away from the injection point.
We can expect that if the temperature is just below $T_c$, the injected fluid will undergo a phase transformation in the high pressure region close to the center of injection.
Therefore, it is possible to find the range of parameters such as injection flux and injection temperature to avoid the phase transformation.
We highlight that we do not consider heat transfer effects in any of these calculations.

\begin{figure}[htb!]
	\begin{center}
		\includegraphics[width=0.65\textwidth]{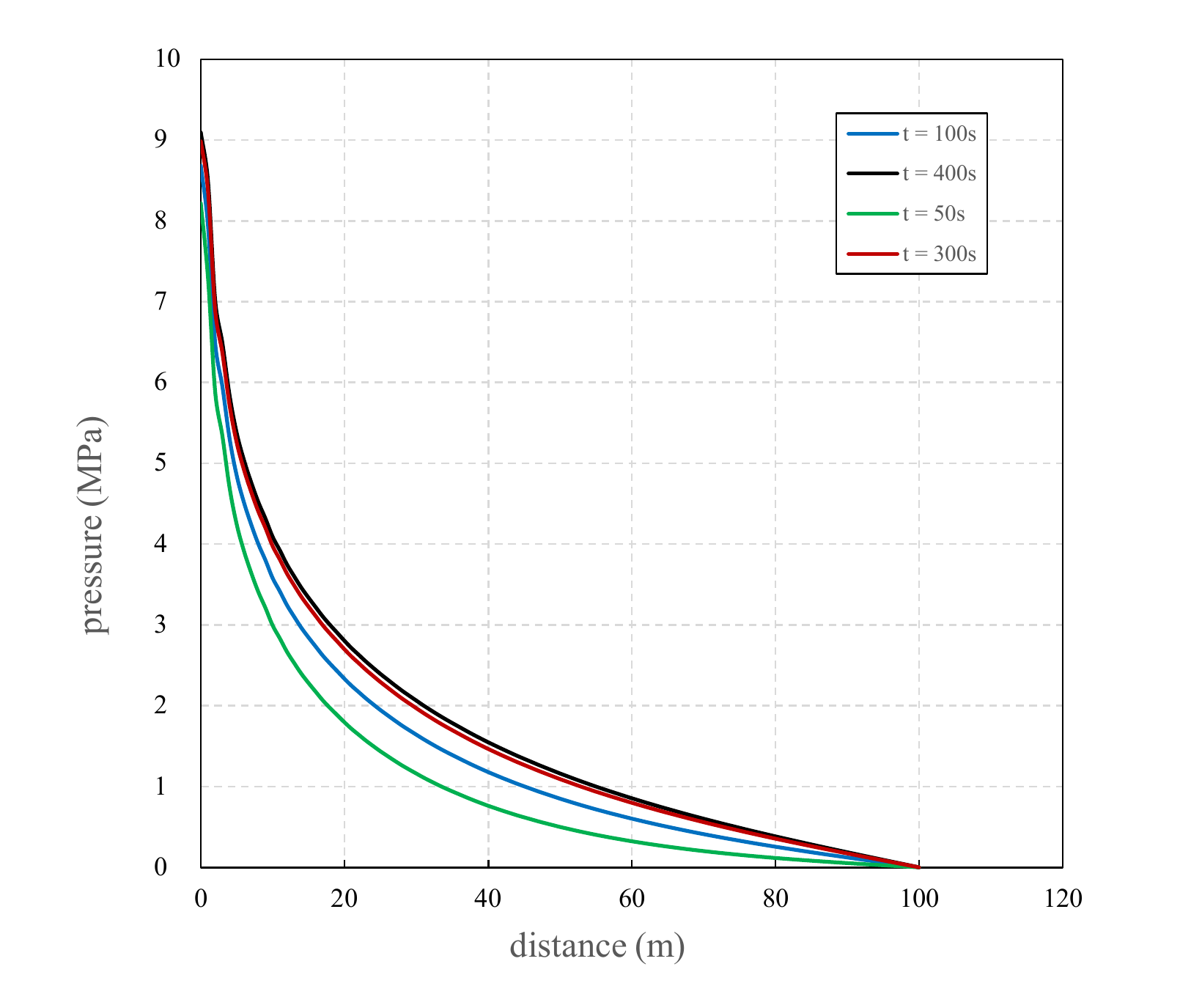}
		\caption{Pore pressure vs. distance at different times.}
		\label{fig:14}
	\end{center}
\end{figure}

To further examine the effect of the pore pressure distribution around the injection point, the predicted pore pressure value at different times is plotted in Fig \ref{fig:14}. 
This simulation uses a flux $\SI{2.2}{\kilo\gram \per \meter^2.\second} $ $\left( \SI{1.112e3}{\liter \per \meter^2.\second} \right)$  and temperature $T = \SI{320}{\kelvin}$ at the injection point. 
The plot clearly shows the increase of pore pressure with time.

\section{Multicomponent Fluid Phase} \label{multiple compressible fluid}

In this section, we apply the model to a porous medium with two immiscible ideal gases, motivated by experimental approaches such as \cite{sarem1966three}.

The Helmholtz free energy density per unit current volume of an ideal gas is:
\begin{equation}
    W_i(\rho_i) = -\rho_i R T \left( 1+\frac{3}{2}\log{\left(\frac{3}{2}R T\right)} - \log{(\rho_i\xi_i)} \right)
\end{equation}
where $R$ and $\xi_i$ are constants and $T$ is the temperature ~\cite{muller2009fundamentals}. 
Following \eqref{eqn:fluid-energy}, we rewrite the energy density of the gases as:
\begin{equation}
    W_{0i}(J,P_{0i}, \phi_i) = 
    \frac{J \phi}{\phi_{0i}} W(J^{-1} P_{0i} / \phi_i) = 
    -\frac{P_{0i}}{\phi_i} R T \left( 1+\frac{3}{2}\log\left(\frac{3}{2}R T\right) - \log\left(\frac{P_{0i}\xi_i}{J\phi_i}\right) \right)
\end{equation}
The strain energy density of the solid phase is modeled as neo-Hookean as in \eqref{eq:Neo-Hookean}.

\begin{figure}[htb!]
	\begin{center}
		\includegraphics[width=0.75\textwidth]{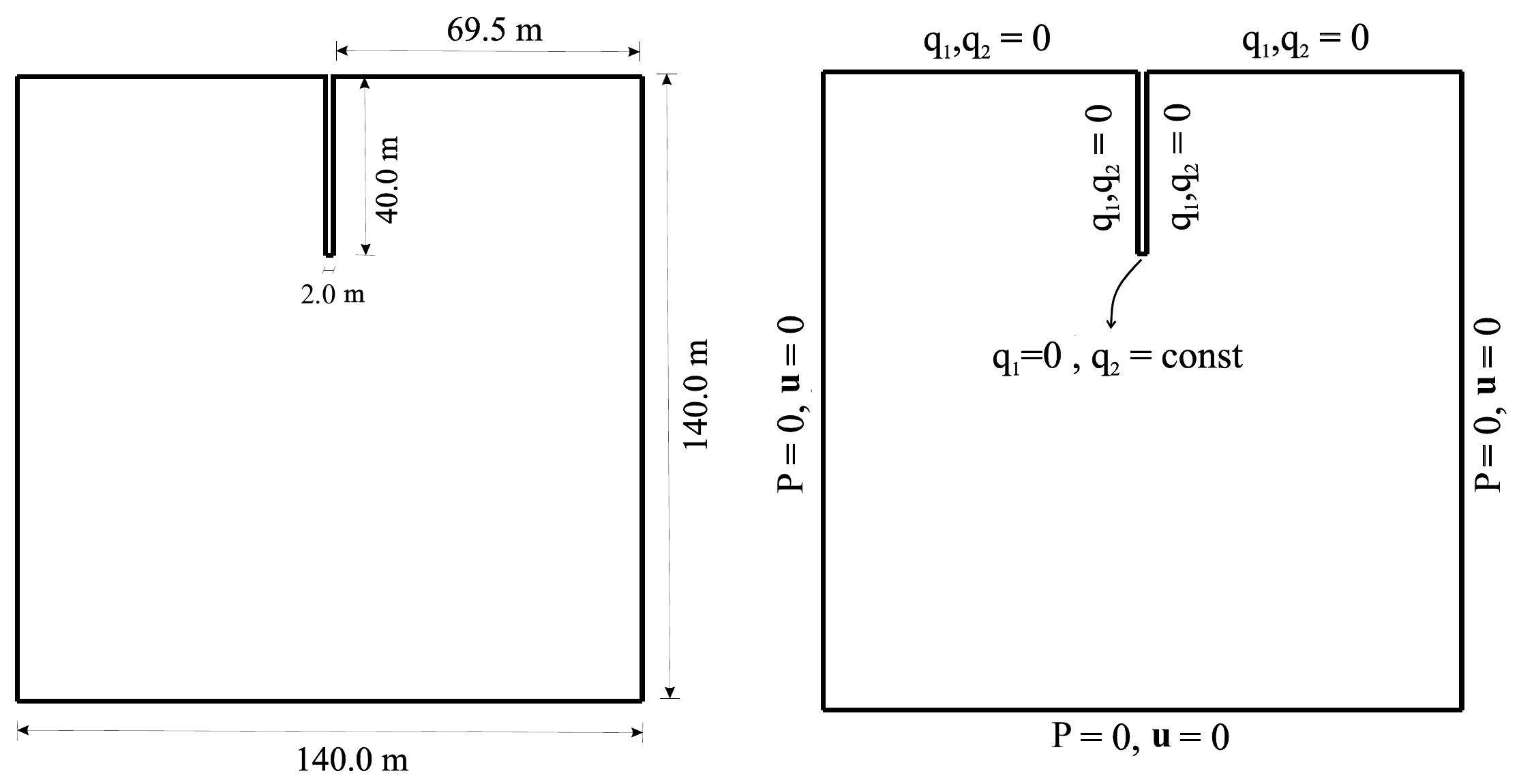}
		\caption{Geometry and boundary conditions.}
		\label{fig:15}
	\end{center}
\end{figure}

By applying \eqref{eq:chemical potential} and \eqref{eq:fluid flux}, we find the chemical potential and flux of the fluid phases:
\begin{align}
    \eta_{0i} & = - R T \left( \frac{3}{2}\log\left(\frac{3}{2}R T\right) - \log\left(\rho_i\xi_i\right) \right)
    \\
    \bfq_i & = -\phi_i\bfk_i R T \nabla \rho_i
\end{align}
The volume fraction of the solid phase is fixed.
Consequently, the unknown variables of system are $\bfF$, $P_{i0}$ and $\phi_i$. 
Therefore, the balance  of mass for the fluid phases and the relation between fluid densities can be found using \eqref{eqn:pressure_equality}:
\begin{align}
    \frac{\dm}{\dm t} P_i J  &= J \divergence \left(\bfk_i R T\phi_i \nabla \rho_i \right)
    \\
    \frac{P_1}{\phi_1} &= \frac{P_2}{\phi_2}
\end{align}
We assume a porous medium with porosity $\phi_{0f} = \num{0.2}$, and relation between the viscosity of two fluids is $\gamma_1 = 2 \gamma_2$, and the geometry and boundary conditions are shown in Fig. \ref{fig:15}.
The properties of fluids and soil are given in Table \ref{tab:3}. 

\begin{figure}[htb!]
	\begin{center}
		\includegraphics[width=0.9\textwidth]{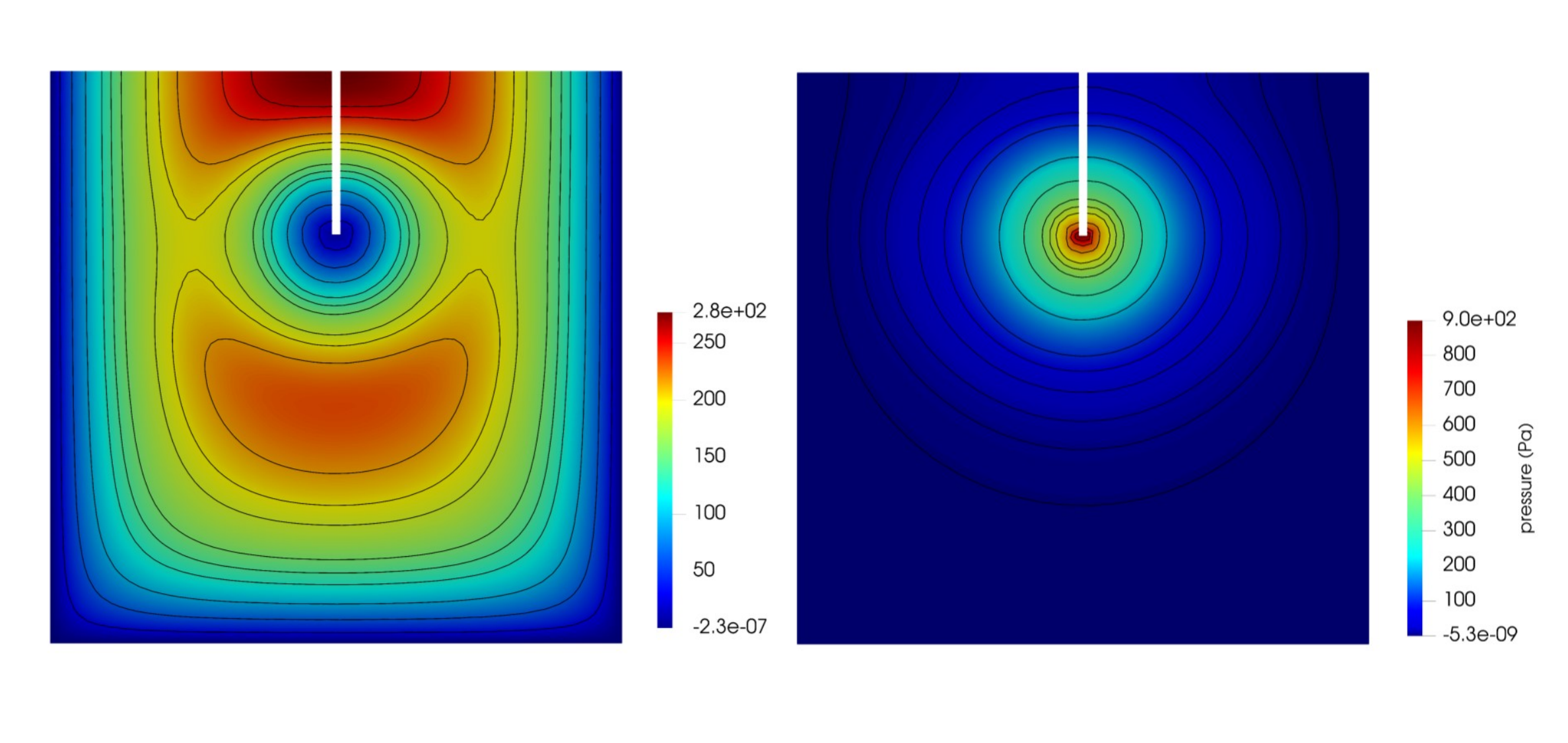}
		\caption{Pore pressure distribution of gas $1$ (left) and gas $2$ (right) at $t=\SI{550}{\second}$.}
		\label{fig:17}
	\end{center}
\end{figure}

At the beginning of the injection ($t = 0$), the domain is saturated with gas $1$; we then begin injecting gas $2$ with the flux $q_2 = \SI{.44}{\kilo\gram \per \meter^2.\second} $ $\left(\SI{2.225e2}{\liter \per \meter^2.\second} \right)$ at the injection point.
Fig. \ref{fig:17} shows the pore pressure distribution of both gases around the injection point at $t=\SI{550}{\second}$. 
Fig. \ref{fig:17} (left) shows that the injection of gas $2$ into the porous medium causes the density of gas $1$ ($P_1$) to decrease around the injection point, and therefore the pore pressure of gas $1$ decreases in this region. 
Further, Fig. \ref{fig:17} (right) shows the increase of the pore pressure due to gas $2$ in regions close to the injection point.

\begin{table*}[h]
    \centering
    \begin{tabular}{l l}
        \hline
        Property & Value\\
        \hline
        Solid phase Lame modulus, $\lambda$ \qquad \qquad & \SI{57.7}{\mega\pascal} \\
        Solid phase Lame modulus, $\mu$ & \SI{38.6}{\mega\pascal} \\
        Hydraulic conductivity, $\tilde{k}_1$  & \SI{0.1}{m\per s}\\
        Hydraulic conductivity, $\tilde{k}_2$  & \SI{0.2}{m\per s}\\
        Gas constant, $R$  &  \SI{8.32}{{\cubic\metre.\pascal}\per {\kelvin.\mole}} \\
        \hline
    \end{tabular}
    \caption{Properties of solid and fluid phases.}\label{tab:3}
\end{table*}

\begin{figure}[htb!]
	\begin{center}
		\includegraphics[width=0.65\textwidth]{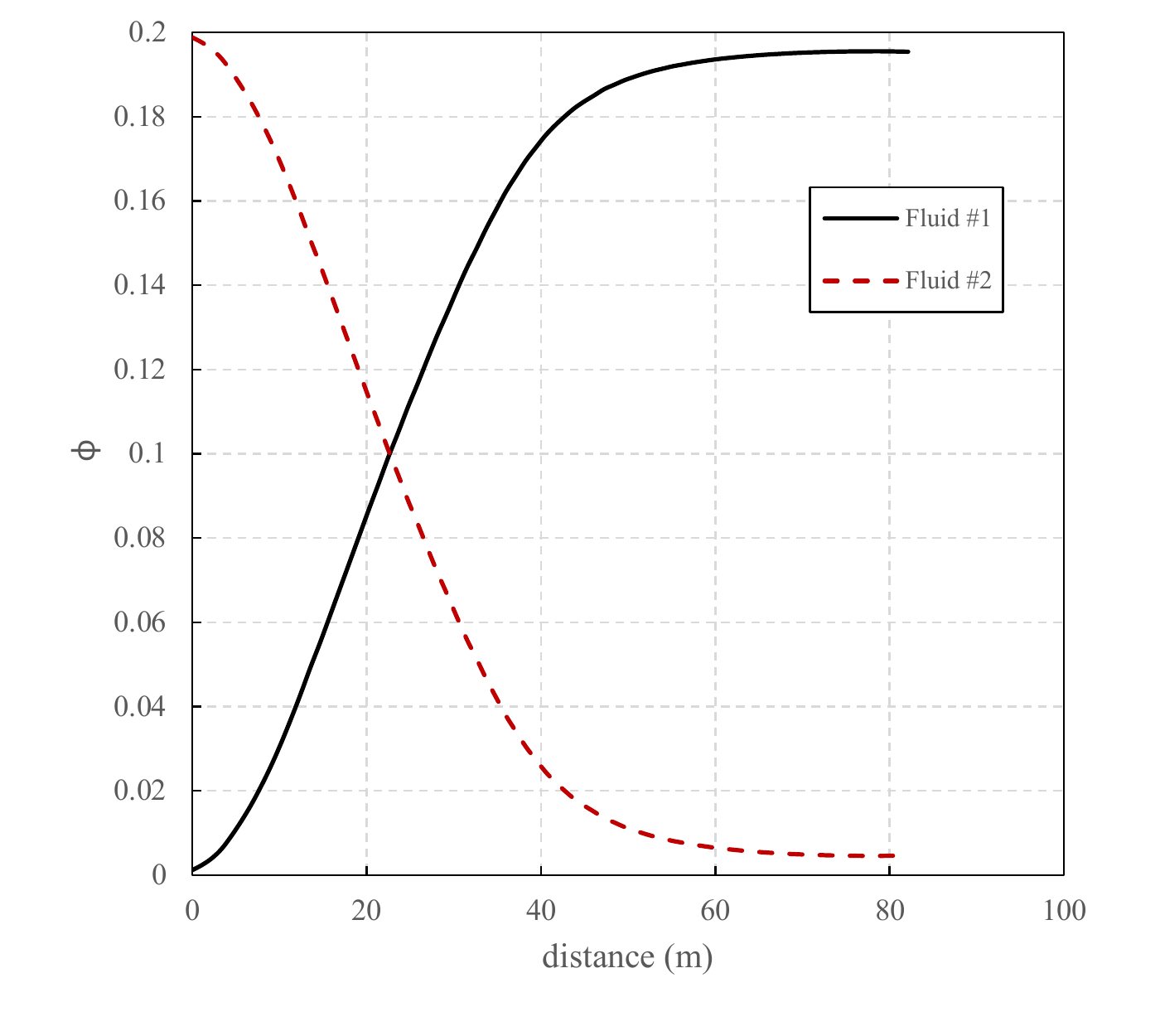}
		\caption{Volume fraction of each fluid component vs. distance from the injection point at $t=\SI{550}{\second}$.}
		\label{fig:18}
	\end{center}
\end{figure}

Due to the compressibility of fluids, the volume fraction of each gas changes with time and distance from the injection point. 
Fig. \ref{fig:18} shows the variation of volume fraction of each gas with vertical location below the injection point, at time ($t=\SI{550}{\second}$). 
After $\SI{550}{\second}$ from the start of the injection, the volume fraction of gas $1$ close to the injection point decreases to $\num{0.001}$ and the volume fraction of gas $2$ increases to $\num{0.199}$.
Since the medium initially is saturated with gas $1$, we find that the volume fraction of gas $1$ tends to $0.2$ far from the injection point.

\section{Unsaturated Medium with an Incompressible Fluid} \label{unsaturated system with incompressible fluid}

In this section, we discuss an unsaturated porous medium consisting of a compressible solid phase, incompressible fluid (liquid) phase, and a compressible fluid (gas) phase, motivated by experimental approaches such as \cite{dria1993three}.

Using the ideal gas model for the compressible fluid phase, the total energy density can be written as:
\begin{equation}
    W_0(\bfF,P_{0g},P_{0l})
    =
    \phi_s W_{0s}(\bfF) 
    - P_{0g} RT\left( 1+\frac{3}{2}\log\left(\frac{3}{2}RT\right) - \log\left(\frac{P_{0g}\xi}{J(1-\phi_s)-\frac{P_{0l}}{\rho_l}}\right) \right)
\end{equation}
Using \eqref{eq:chemical potential}, the chemical potential of gas and liquid can be derived as:
\begin{align}
    \eta_{0g} & = -RT \left( \frac{3}{2}\log\left(\frac{3}{2}RT\right)-\log \left(\frac{\rho_{0g}\xi}{J(1-\phi_s)-\frac{P_{0l}}{\rho_l}} \right) \right)
    \\
    \eta_{0l} & = \frac{\rho_g RT}{\rho_l} = \frac{P_g RT}{(1-\phi_s-\phi_l)\rho_l }
\end{align}
Further, the flux of gas and liquid can be written as:
\begin{align}
    \bfq_g &= -\phi_g \bfk_g RT \nabla\rho_g
    \\
    \bfq_l &= -\phi_l \bfk_l RT \nabla\rho_l
\end{align}
The geometry is the same as Section \ref{multiple compressible fluid}, and the properties of the solid and gas phases are listed in Table \ref{tab:4}.
We consider a porous medium with initial solid volume fraction $\phi_{0s} = 0.9$, saturated with the ideal gas ($\phi_{0g} = 0.1$ and $p_g=\SI{3300}{\pascal}$ at $t = 0$).
We then inject the liquid with the flux $q_l = \SI{200}{\liter \per \meter^2.\second}$. 

\begin{table*}[h]
    \centering
    \begin{tabular}{l l}
        \hline
        Property & Value\\
        \hline
        Solid phase Lame constant, $\lambda$ \qquad \qquad &  \SI{57.7}{\mega\pascal}\\ 
        Solid phase Lame constant, $\mu$ & \SI{38.6}{\mega\pascal}\\
        True permeability, $\kappa$ & \num{1.8e-7} \si[per-mode=symbol]{\meter^2}\\
        Gas constant, $R$ & \SI{8.32}{{\cubic\metre.\pascal}\per {\kelvin.\mole}} \\
        Gas viscosity, $\gamma_g$ & \num{1.8e-5} \si[per-mode=symbol]{\pascal.\second}\\
        Liquid density, $\rho_l$ & \num{1000} \si[per-mode=symbol]{\kilo\gram\per\cubic\meter}\\
        Liquid viscosity, $\gamma_l$ & \SI{0.001}{\pascal.\second}\\
        \hline
    \end{tabular}
    \caption{Properties of the solid and fluid phases.}\label{tab:4}
\end{table*}

\begin{figure}[htb!]
	\begin{center}
		\includegraphics[width=0.9\textwidth]{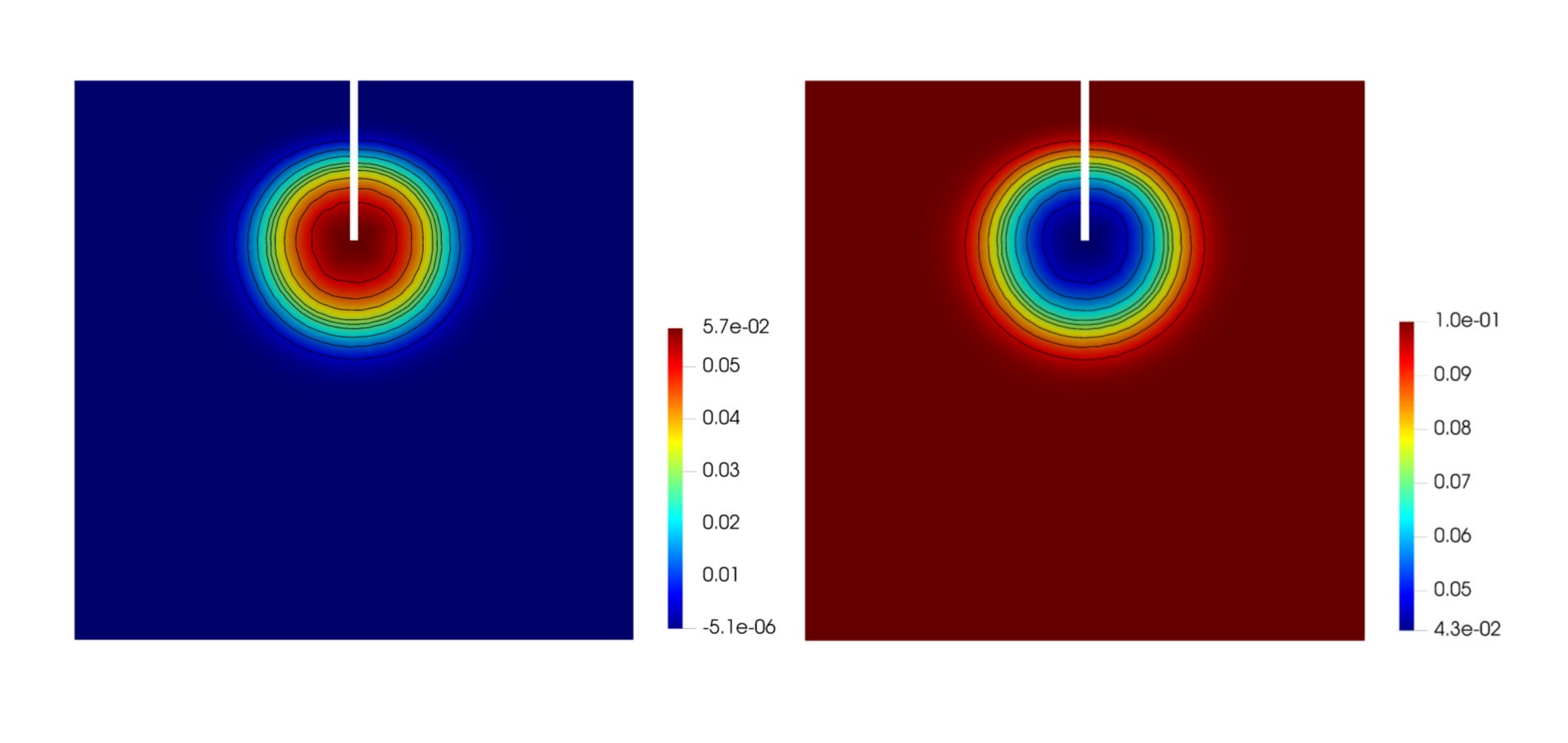}
		\caption{Volume fraction of liquid (left) and gas (right) around the injection point at $t=\SI{200}{\second}$.}
		\label{fig:20}
	\end{center}
\end{figure}

\begin{figure}[htb!]
	\begin{center}
		\includegraphics[width=0.55\textwidth]{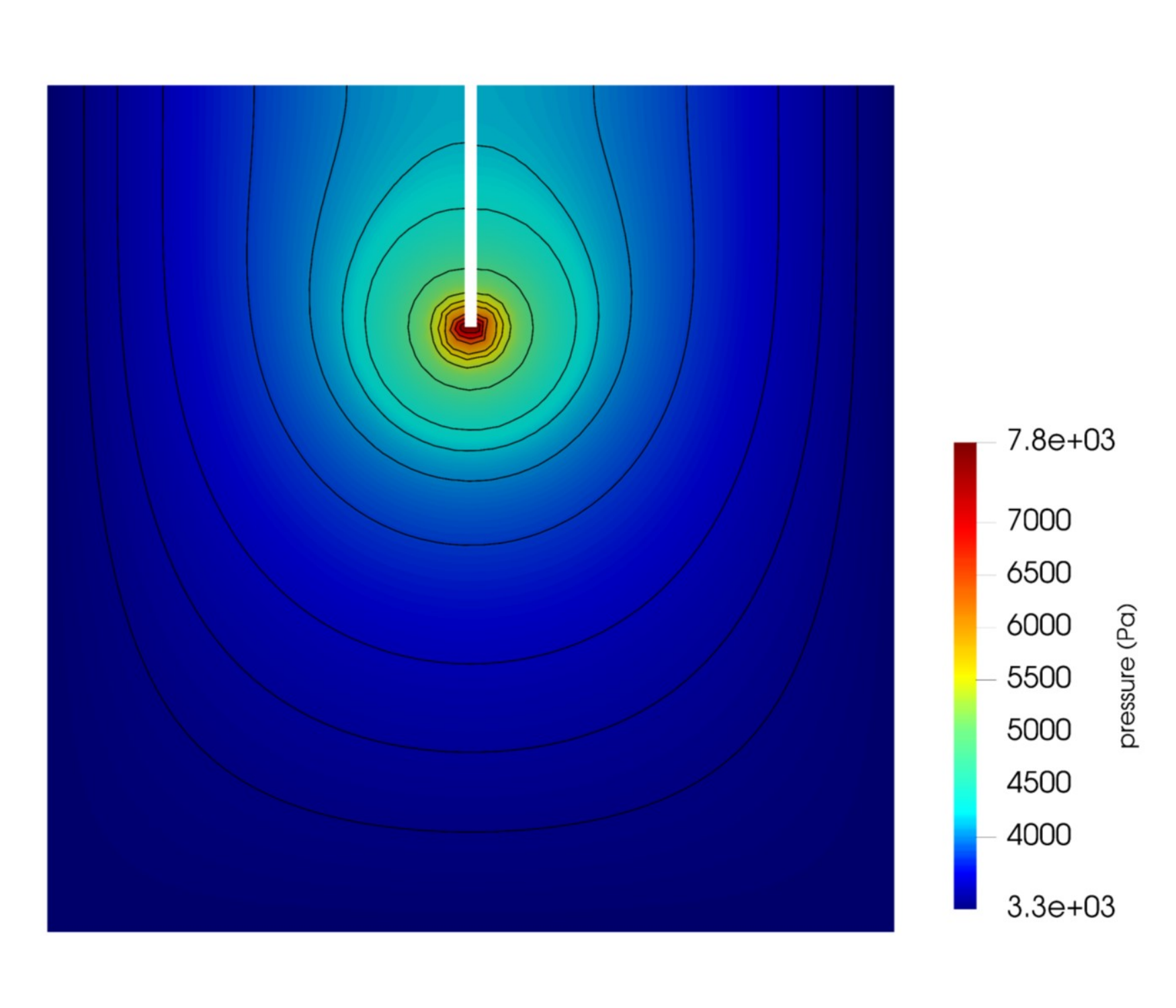}
		\caption{The pore pressure distribution around the injection point at $t=\SI{200}{\second}$.}
		\label{fig:21}
	\end{center}
\end{figure}

Fig. \ref{fig:20} shows the volume fraction distributions for the gas and liquid phases.
As expected, injecting the liquid into the porous medium with gas causes the liquid to displace the gas in the pore spaces, and the volume fraction of gas decreases in the vicinity of the injection point.
Also, the pore pressure distribution due to the liquid injection is shown in Fig. \ref{fig:21}; we see that the highest value of pore pressure is in the vicinity of the injection point, and tends gradually to the far-field pressure of the gas.

\section{Concluding Remarks} \label{Concluding Remarks}

In this study, we have proposed a variational energy-based continuum mechanics framework in the large deformation regime to model coupled fluid transport and deformation in a porous deformable medium. 
The free energy density function of the porous medium is defined as additively composed of the free-energies of the solid and fluid phases of the medium.
The variational structure provides for the standard use of the Finite Element Method to solve the proposed model; in general, variational approaches can provide important advantages for numerical solution \cite{oden2012variational}.

The proposed formulation was tested on two benchmark problems, namely Terzaghi's 1-d and Mandel's 2-d consolidation problems.
Further, the model enabled the modeling of porous media containing a single (gas-liquid phase transforming) and multiple immiscible fluid phases.
Furthermore, the method was used to model the behavior of an unsaturated porous medium.
As a specific application of the proposed model, we investigated fluid injection into soil.

An important advantage of the variational formulation, that we are currently investigating and will be reported in the future, is in enabling an overall variational structure to the problem of combining a physics-based approach with a data-driven approach through a Bayesian framework.
In addition, plastic effects are important in many realistic situations; it would therefore be useful to extend this model to incorporate powerful variational models for plasticity that have been proposed, e.g. \cite{weinberg2006variational,ortiz2004variational}.


\section*{Author Contributions}

All authors contributed to the formulation of the model.
Karimi, Walkington, and Dayal contributed to the numerical implementation.
All authors contributed to the writing of the paper.

\begin{acknowledgments}
    We thank the National Science Foundation for support through XSEDE resources provided by Pittsburgh Supercomputing Center.
    Mina Karimi acknowledges financial support from the Scott Institute.
    Noel Walkington acknowledges financial support from NSF (DMS 1729478 and DMS 2012259).
    Kaushik Dayal acknowledges financial support from NSF (CMMI MOMS 1635407, DMREF 1628994), ARO (MURI W911NF-19-1-0245), ONR (N00014-18-1-2528), BSF (2018183), and an appointment to the National Energy Technology Laboratory sponsored by the U.S. Department of Energy.
    This work was funded (in part) by the Dowd Fellowship from the College of Engineering at Carnegie Mellon University; we thank Philip and Marsha Dowd for their financial support and encouragement. 
    We thank Basant Sharma for alerting to us to typos in a previous version.
\end{acknowledgments}

\section*{Data Availability Statement}

A version of the code developed for this work is available at \url{github.com/minakari/Poromechanics}.



\bibliographystyle{alpha}
\bibliography{poromech-refs}

\end{document}